\pdfoutput=1 
\documentclass[graybox, dvipsnames]{svmult}
\usepackage{mathptmx}       
\usepackage{helvet}         
\usepackage{courier}        
\usepackage{type1cm}        
\usepackage{makeidx}         
\usepackage{graphicx}        
\usepackage{multicol}        
\usepackage[bottom]{footmisc}
\makeindex            
\usepackage{amssymb}
\usepackage{bbm}
\usepackage{algorithm}
\usepackage{algpseudocode}
\usepackage{amssymb}
\usepackage{graphicx}
\usepackage{color}
\usepackage{mathtools}
\usepackage{amsmath}
\usepackage{subfigure}
\usepackage{subfigmat}
\usepackage{amsfonts}
\usepackage{xcolor}
\usepackage{color,soul}
\usepackage{wasysym}
\usepackage[colorlinks, citecolor=Red,linkcolor=Green, breaklinks=true]{hyperref}
\usepackage{bm}

\graphicspath{{./}{./figs/}}

  \def\clap#1{\hbox to 0pt{\hss#1\hss}}

\providecommand{\mat}[1]{\bm{#1}}%
\renewcommand{\vec}[1]{\mathbf{#1}}

\providecommand{\mA}{\ensuremath{\mat{A}}}

\providecommand{\mG}{\ensuremath{\mat{G}}}

\providecommand{\mP}{\ensuremath{\mat{P}}}
\providecommand{\mQ}{\ensuremath{\mat{Q}}}
\providecommand{\mR}{\ensuremath{\mat{R}}}

\providecommand{\va}{\ensuremath{\vec{a}}}
\providecommand{\vb}{\ensuremath{\vec{b}}}

\providecommand{\ve}{\ensuremath{\vec{e}}}
\providecommand{\vf}{\ensuremath{\vec{f}}}

\providecommand{\vp}{\ensuremath{\vec{p}}}

\providecommand{\vv}{\ensuremath{\vec{v}}}
\providecommand{\vw}{\ensuremath{\vec{w}}}
\providecommand{\vx}{\ensuremath{\vec{x}}}

\providecommand{\vz}{\ensuremath{\vec{z}}}

\begin{document}

\title*{Quadrature Strategies for Constructing Polynomial Approximations}
\author{Pranay Seshadri, Gianluca Iaccarino, Tiziano Ghisu}
\institute{Pranay Seshadri \at Department of Engineering, University of Cambridge, Cambridge, U. K., \email{ps583@cam.ac.uk}
\and Gianluca Iaccarino \at Department of Mechanical Engineering and Institute for Computational and Mathematical Engineering, Stanford University, Stanford, California, U. S. A., \email{jops@stanford.edu}
\and Tiziano Ghisu \at Department of Mechanical, Chemical and Materials Engineering, Universit\'{a} di Cagliari, Cagliari, Sardinia, Italy, \email{t.ghisu@unica.it}}
%
%
\maketitle

\abstract{Finding suitable points for multivariate polynomial interpolation and approximation is a challenging task. Yet, despite this challenge, there has been tremendous research dedicated to this singular cause. In this paper, we begin by reviewing classical methods for finding suitable \emph{quadrature} points for polynomial approximation in both the univariate and multivariate setting. Then, we categorize recent advances into those that propose a new sampling approach and those centered on an optimization strategy. The sampling approaches yield a favorable discretization of the domain, while the optimization methods pick a subset of the discretized samples that minimize certain objectives. While not all strategies follow this two-stage approach, most do. Sampling techniques covered include subsampling quadratures, Christoffel, induced and Monte Carlo methods. Optimization methods discussed range from linear programming ideas and Newton's method to greedy procedures from numerical linear algebra. Our exposition is aided by examples that implement some of the aforementioned strategies.}

\section{Introduction}
\label{sec:1}
Over the past few years there has been a renewed research effort aimed at constructing stable polynomial approximations: understanding their theoretical properties \cite{shin2016near, zhou2015weighted, shin2016nonadaptive, seshadri2017effectively, chkifa2015discrete, migliorati2015analysis, cohen2013stability, hampton2015coherence, adcock2018approximating} and extending their practical use. The latter spans from applications in spectral methods \cite{trefethen2000spectral}, uncertainty quantification and the related sensitivity analysis \cite{pettersson2015polynomial, seshadri2014leakage, ghisu2017asme, gmd-2017-107} to dimension reduction \cite{seshadri2018turbomachinery} and design under uncertainty \cite{seshadri2016density}. One topic germane to the field has been to find suitable points for multivariate polynomials using as few points as possible. The application-centric high level idea can be described as follows:
\begin{align*}
\begin{array}{c}
\text{select a polynomial basis}\\
\Downarrow\\
\text{evaluate model at quadrature points }\\
\Downarrow\\
\textrm{estimate polynomial coefficients} \\
\Downarrow\\
\textrm{compute moments, sensitivities and probability density functions} 
\end{array}
\end{align*}

Motivated by reducing the number of samples, and providing scalable computational methods for doing so, numerous deterministic and randomized sampling schemes have been proposed. Before we delve into recent strategies for estimating these using perspectives from least squares, a brief overview of some of the more fundamental concepts is in order. 

Quadrature rules given by $\left\{ \left(\bm{\zeta}_{i},\omega_{i}\right)\right\} _{i=1}^{m}$ seek to find $d-$dimensional points $\bm{\zeta}_{i} \in \mathbb{R}^{d}$ and weights $\omega_{i} > 0$ in $\mathbb{R}$ such that the integral of a function $f$ may be expressed as
\begin{equation}
\int_{\mathbb{R}^{d}} f \left( \bm{\zeta} \right) \rho\left( \bm{\zeta} \right)d \bm{\zeta} \approx\sum_{i=1}^{m}f\left(\bm{\zeta}_{i}\right)\omega_{i}
\label{equ1}
\end{equation}
where $\rho(\bm{\zeta})$ is a known multivariate weight function---i.e., Gaussian, uniform, Cauchy, etc. or product thereof. These $d-$dimensional points are sampled over the support of $\bm{\zeta}=(\zeta^{(1)}, \ldots, \zeta^{(d)} )$, mutually independent random variables under the joint density $\rho(\bm{\zeta})$. By construction, this quadrature rule must also satisfy
\begin{equation}
\sum_{i=1}^{m} \bm{\psi}_{\bm{p} }\left( \bm{\zeta}_{i}\right) \bm{\psi}_{\bm{q} }\left( \bm{\zeta}_{i}\right)\omega_{i}\approx\int_{\mathbb{R}^{d}} \bm{\psi}_{\bm{p}}\left( \bm{\zeta}  \right) \bm{\psi}_{\bm{q}}\left(\bm{\zeta} \right)\rho\left(\bm{\zeta} \right)d\bm{\zeta}=\delta_{\bm{pq}}
\label{equ2}
\end{equation}
where $\bm{\psi}_{\bm{p}} (\bm{\zeta})$ is a multivariate polynomial $L^{2}$-orthogonal on $\mathbb{R}^{d}$ when weighted by the joint density $\bm{\rho}(\bm{\zeta})$. Here $\delta_{\bm{pq}}$ denotes the Kronecker delta; subscripts $\bm{p}$ and $\bm{q}$ are multi-indices that denote the order of $\bm{\psi}$ and its composite univariate polynomials $\psi_{j}$ via
\begin{equation}
\bm{\psi}_{\bm{p}}\left( \bm{\zeta} \right)=\prod_{k=1}^{d}\psi_{p_{k}}^{\left(k\right)}\left(\zeta^{\left(k\right)}\right) \; \; \text{where} \;  \; \bm{p}=\left(p_{1},\ldots,p_{d}\right)\in\mathbb{N}^{d}.
\end{equation}
The finite set of indices in $\bm{p}$ is said to belong to a multi-index set $\mathcal{I}$. The number of multi-indices present in $\bm{p}, \bm{q} \in \mathcal{I}$ is set by choosing elements in $\mathcal{I}$ to follow certain rules, e.g., the sum of all univariate polynomial orders must satisfy $\sum_{i=1}^{d}p_{i}\leq k$, yielding a \emph{total order} index set \cite{xiu2010numerical}. A \emph{tensor order} index, which scales exponentially in $d$ is governed by the rule $max_{k}p_{k}\leq k$. Other well-known multi-index sets include hyperbolic cross \cite{zhou2015weighted} and hyperbolic index sets \cite{blatman2011adaptive}. We shall denote the number of elements in each index set by $n=\text{card} \left( \bm{p} \right)$.

From \eqref{equ1} and \eqref{equ2} it should be clear that one would ideally like to minimize $m$ and yet achieve equality in the two expressions. Furthermore, the \emph{degree of exactness} associated with evaluating the integral of the function in \eqref{equ1} will depend on the highest order polynomial (in each of the $d-$dimensions) that yields equality in \eqref{equ2}. Prior to elaborating further on multivariate quadrature rules, it will be useful to detail key ideas that underpin univariate quadrature rules. It is important to note that there is an intimate relationship between polynomials and quadrature rules; ideas that date back to Gauss, Christoffel and Jacobi.  

\subsection{Classical quadrature techniques}
Much of the development of classical quadrature techniques is due to the foundational work of Gauss, Christoffel and Jacobi. In 1814 Gauss originally developed his quadrature formulas, leveraging his theory of continued fractions in conjunction with hypergeometric series using what is today known as Legendre polynomials. Jacobi's contribution was the connection to \emph{orthogonality}, while Christoffel is credited with the generalization of quadrature rules to non-negative, integrable weight functions $\rho(\zeta) \geq 0$ \cite{gautschi1981survey}. The term \emph{Gauss-Christoffel} quadrature broadly refers to all rules of the form \eqref{equ1} that have a degree of exactness of $(2m-1)$ for $m$ points when $d=1$. The points and weights of Gauss rules may be computed either from the moments of a weight function or via the three-term recurrence formula associated with the orthogonal polynomials. In most codes today, the latter approach, i.e., the Golub and Welch \cite{golub1969calculation} approach is adopted. It involves computing the eigenvalues of a tridiagonal matrix\footnote{Known colloquially as the Jacobi matrix.}---incurring complexity $\mathcal{O}\left( m^2 \right)$---of the recurrence coefficients associated with a given orthogonal polynomial. For uniform weight functions $\rho(\zeta)$ these recurrence coefficients are associated with Legendre polynomials, and the resulting quadrature rule is termed Gauss-Legendre. When using Hermite polynomials---orthogonal with respect to the Gaussian distribution---the resulting quadrature is known as Gauss-Hermite. In applications where the weight functions are arbitrary, or data-driven\footnote{In the case of data-driven distributions, kernel density estimation or even a maximum likelihood estimation may be required to obtain a probability distribution that can be used by the discretized Stieltjes procedure.}, the discretized Stieltjes procedure (see section 5 in \cite{gautschi1985orthogonal}) may be used for computing the recurrence coefficients. 

\begin{figure}
\begin{subfigmatrix}{3}
\subfigure[]{\includegraphics[]{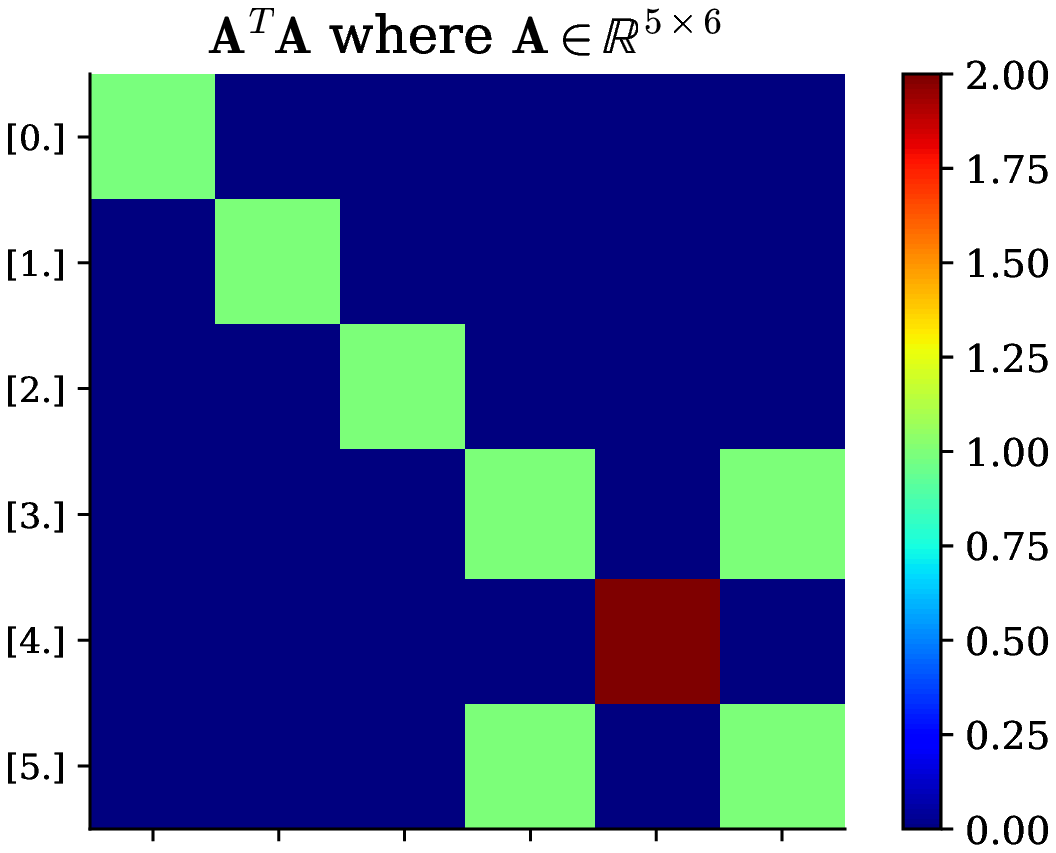}}
\subfigure[]{\includegraphics[]{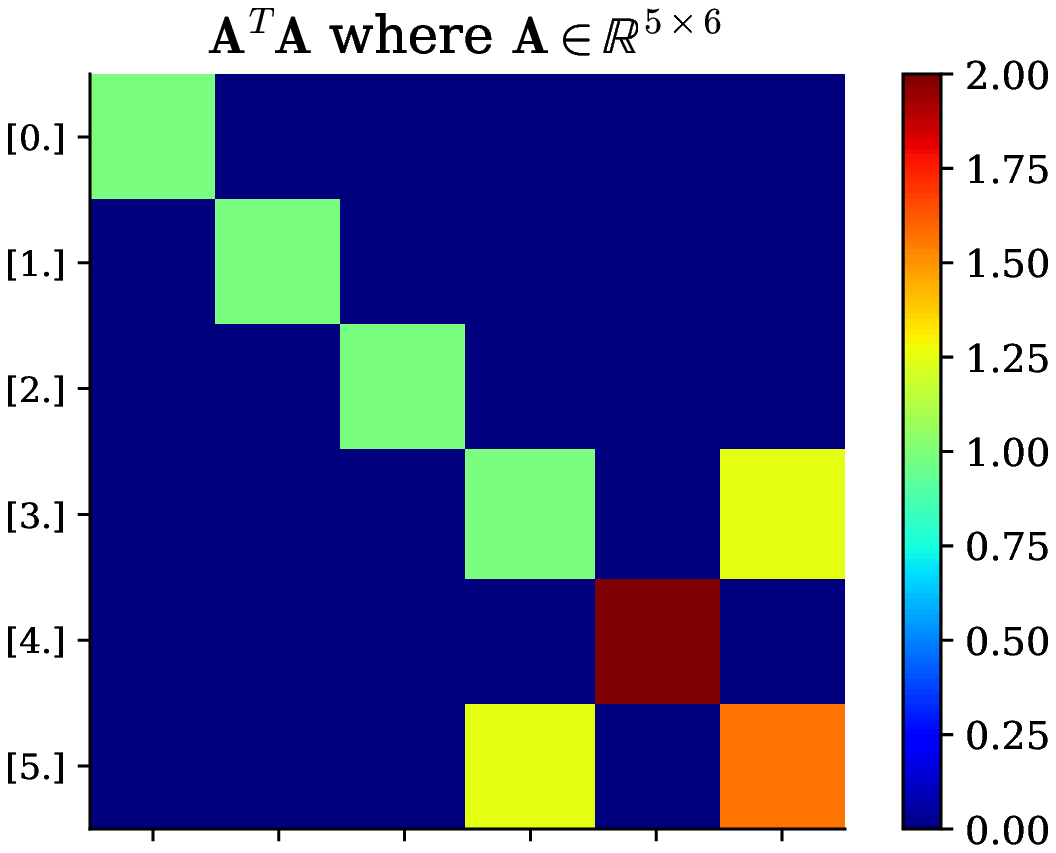}}
\subfigure[]{\includegraphics[]{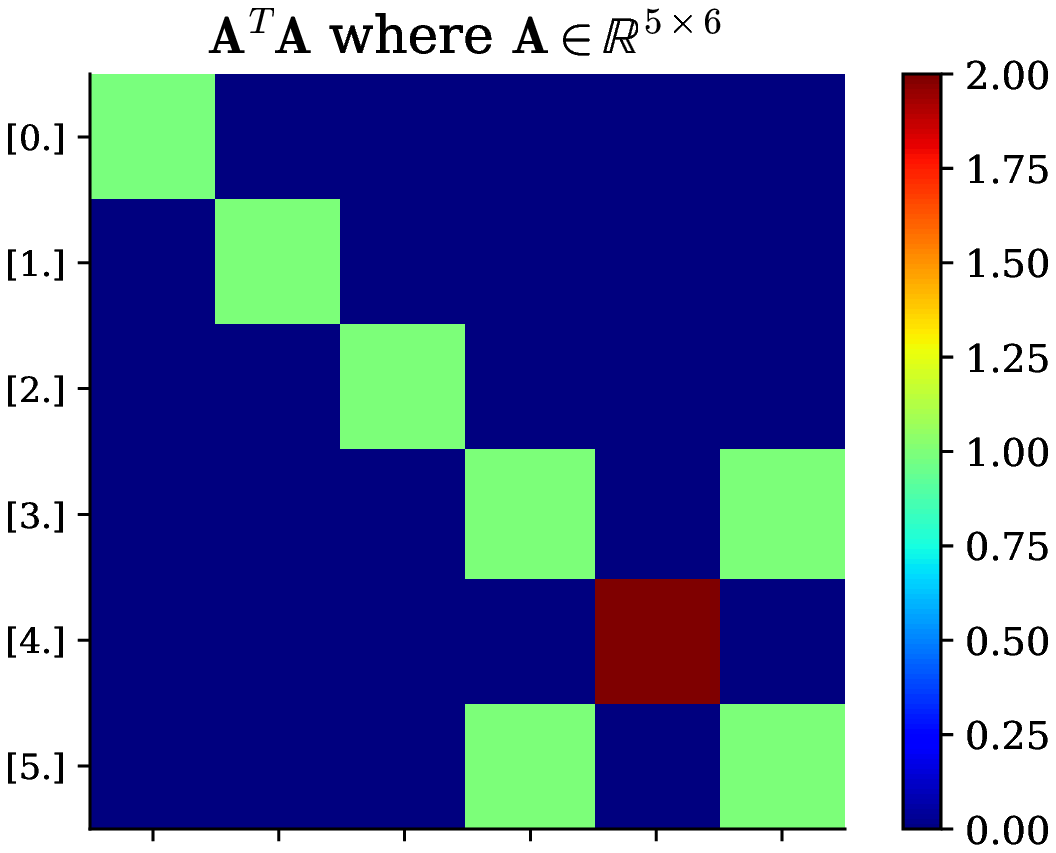}}
\end{subfigmatrix}
\caption{The Gram matrix $\mG$ for three quadrature rules showing their practical degrees of exactness: (a) Gauss-Legendre; (b) Gauss-Lobatto; (c) Clenshaw-Curtis. The vertical axis in the figures show the polynomial degree and the colorbar shows the values of the entries of these matrices.}
\label{doe}
\end{figure}

\begin{svgraybox}
\textbf{Degree of exactness:}
The notion of degree of exactness dates back to Radau \cite{gautschi1981survey}; all univariate quadrature rules are associated with a \emph{degree of exactness}, referring to the highest degree polynomial that can be integrated exactly for a given number of points $m$. For instance, Gauss-Legendre quadrature has a degree of exactness of $(2m-1)$. Gauss-Lobatto quadrature rules on the other hand have a degree of exactness of $(2m-3)$, while Clenshaw-Curtis (see \cite{gentleman1972implementing}) have a degree of exactness of $(m-1)$---although in practice comparable accuracy to Gauss quadrature rules can be obtained (see Trefethen's monograph \cite{trefethen2008gauss} and Davis and Rabinowitz  \cite{davis2007methods}). One way to interpret this degree of exactness is to inspect the elements of a Gram matrix $\mG = \mA^T \mA$, where $\mA$ is formed by evaluating the orthogonal polynomials at the quadrature points
\begin{equation}
\mA(i,j) = \psi_{j} \left(\zeta_{i}\right) \sqrt{\omega}_{i} \; \; \text{where}\;  \mA \in \mathbb{R}^{m \times n},
\label{def_of_A}
\end{equation}
with $m$ quadrature points and the first $n$ polynomials. Thus, each element of the Gram matrix seeks to approximate
\begin{equation}
\mG(p,q) = \int\psi_{p}\left(\zeta\right)\psi_{q}\left(\zeta\right)\rho\left(\zeta\right)d\zeta \approx \sum_{i=1}^{m}\psi_{p}\left(\zeta_{i}\right)\psi_{q}\left(\zeta_{i}\right)\omega_{i}=\delta_{pq}
\label{doe_1}
\end{equation}
To clarify this, consider the example case of a Gauss-Legendre rule with $m=5$. The highest polynomial degree that this rule can integrate up to is $9$, implying that the first 4-by-4 submatrix of $\mG$ will be the identity matrix as the combined polynomial degree of the terms inside the integral in \eqref{doe_1} is 8. This is illustrated in Figure~\ref{doe}(a). In (b) and (c) similar results are shown for Gauss-Lobatto (the highest degree being 7)  and Clenshaw-Curtis (which integrates higher than degree 5). For all the aforementioned subfigures, element-wise deviations from the identity can be interpreted as the \emph{internal aliasing} errors associated with each quadrature rule. 
\end{svgraybox}

Classical extensions of Gauss quadrature include:
\begin{itemize}
\item \textbf{Gauss-Radau and Gauss-Lobatto}: The addition of end-points $(a,b)$ in quadrature rules. In Gauss-Radau the addition of one end-point $(a,b]=\left\{ \zeta\in\rho|a<\zeta\leq b\right\}$ and both in Gauss-Lobatto, i.e., $[a,b]=\left\{ \zeta\in\rho|a\leq\zeta\leq b\right\}$ \cite{gautschi1981survey}.
\item \textbf{Gauss-Kronrod and Gauss-Patterson}: Practical extensions of $m$-point Gauss rules by adaptively adding $m+1$ points to the quadrature rule to yield a \emph{combined} degree of exactness of $(3m+1)$ or $(3n+2)$ depending on whether m is even or odd respectively; based on the work of Kronrod \cite{laurie1997calculation}. Patterson \cite{patterson1968optimum} rules are based on optimizing the degree of exactness of Kronrod points and with extensions to Lobatto rules. The intent behind these rules is that they admit certain favorable \emph{nesting} properties.
\item \textbf{Gauss-Tur\'{a}n}:  The addition of derivative values (in conjunction with function values) within the integrand (see \cite{gautschi2004orthogonal}).
\end{itemize}

\subsection{Multivariate extensions of classical rules}
Multivariate extensions of univariate quadrature rules exist in the form of tensor product grids, assuming one is ready to evaluate a function at $m^{d}$ points. They are expressed as:
\begin{align}
\begin{split}
\int_{\mathbb{R}^{d}}f\left( \bm{\zeta} \right)\rho\left( \bm{\zeta}  \right)d \bm{\zeta} & \approx  \left(\mathcal{Q}_{q}^{m_{1}}\otimes\ldots\otimes\mathcal{Q}_{q}^{m_{d}}\right)\left( f \right) \\
& = \sum_{j_{1}=1}^{m_{1}}\ldots\sum_{j_{d}=1}^{m_{d}}f\left(\zeta_{j_{d}},\ldots,\zeta_{j_{d}}\right)\omega_{j_{1}}\ldots\omega_{j_{d}} \\
& = \sum_{j=1}^{m}f\left(\bm{\zeta}_{j}\right)\omega_{j}.
\end{split}
\end{align}
The notation $\mathcal{Q}_{q}^{m_{i}}$ denotes a linear operation of the univariate quadrature rule applied to $f$ along direction $i$; the subscript $q$ stands for \emph{quadrature}. 

In the broader context of approximation, one is often interested in the projection of $f$ onto $\psi_{\bm{p}} \left( \bm{\zeta} \right)$. This \emph{pseudospectral approximation} is given by
\begin{align}
f\left( \bm{\zeta}  \right)\approx\sum_{i=1}^{n} x_{i} \bm{\psi_{i}}  \left( \bm{\zeta} \right), & & \text{where} & & x_{i}=\sum_{j=1}^{m}f\left(\bm{\zeta}_{j}\right)\bm{\psi}_{i}\left(\bm{\zeta}_{j}\right)\omega_{j}.
\end{align}
Inspecting the decay of these \emph{pseudospectral coefficients}, i.e., $\vx = \left(x_{1},\ldots,x_{m}\right)^{T}$ is useful for analyzing the quality of the overall approximation to $f\left( \bm{\zeta}  \right)$, but more specifically for gaining insight into which directions the function varies the greatest, and to what extent. We remark here that for simulation-based problems where an approximation to $f\left( \bm{\zeta}  \right)$ may be required, it may be unfeasible to evaluate a model at a tensor grid of quadrature points, particularly if $d$ is large. Sparse grids \cite{ gerstner1998numerical, smolyak1963quadrature} offer moderate computational attenuation to this problem. One can think of a sparse grid as linear combinations of \emph{select} anisotropic tensor grids
\begin{equation}
\int_{\mathbb{R}^{d}}f\left( \bm{\zeta} \right)\rho\left( \bm{\zeta}  \right)d \bm{\zeta} \approx \sum_{\bm{r} \in\mathcal{K}} \bm{\alpha}  \left(\bm{r}  \right)\left(\mathcal{Q}_{q}^{m_{1}}\otimes\ldots\otimes\mathcal{Q}_{q}^{m_{d}}\right) \left( f \right) ,
\end{equation}
where for a given \emph{level} $l$---a variable that controls the density of points and the highest order of the polynomial in each direction---the multi-index $\mathcal{K}$ and coefficients $\bm{\alpha}(\bm{r})$ are given by
\begin{align}
\mathcal{K}=\left\{ \bm{r}  \in\mathbb{N}^{d}:l+1\leq\left| \bm{r}  \right|\leq l+d\right\} \; \;  \text{and} \; \;  \bm{\alpha}\left( \bm{r} \right)=\left(-1\right)^{  -\left|\bm{r} \right|+d+l  } \binom{d-1}{-\left|\bm{r} \right|+d+l }.
\end{align}
Formulas for the pseudospectral approximations via sparse grid integration can then be written down; we omit these in this exposition for brevity.

\begin{figure}
\begin{subfigmatrix}{2}
\subfigure[]{\includegraphics[]{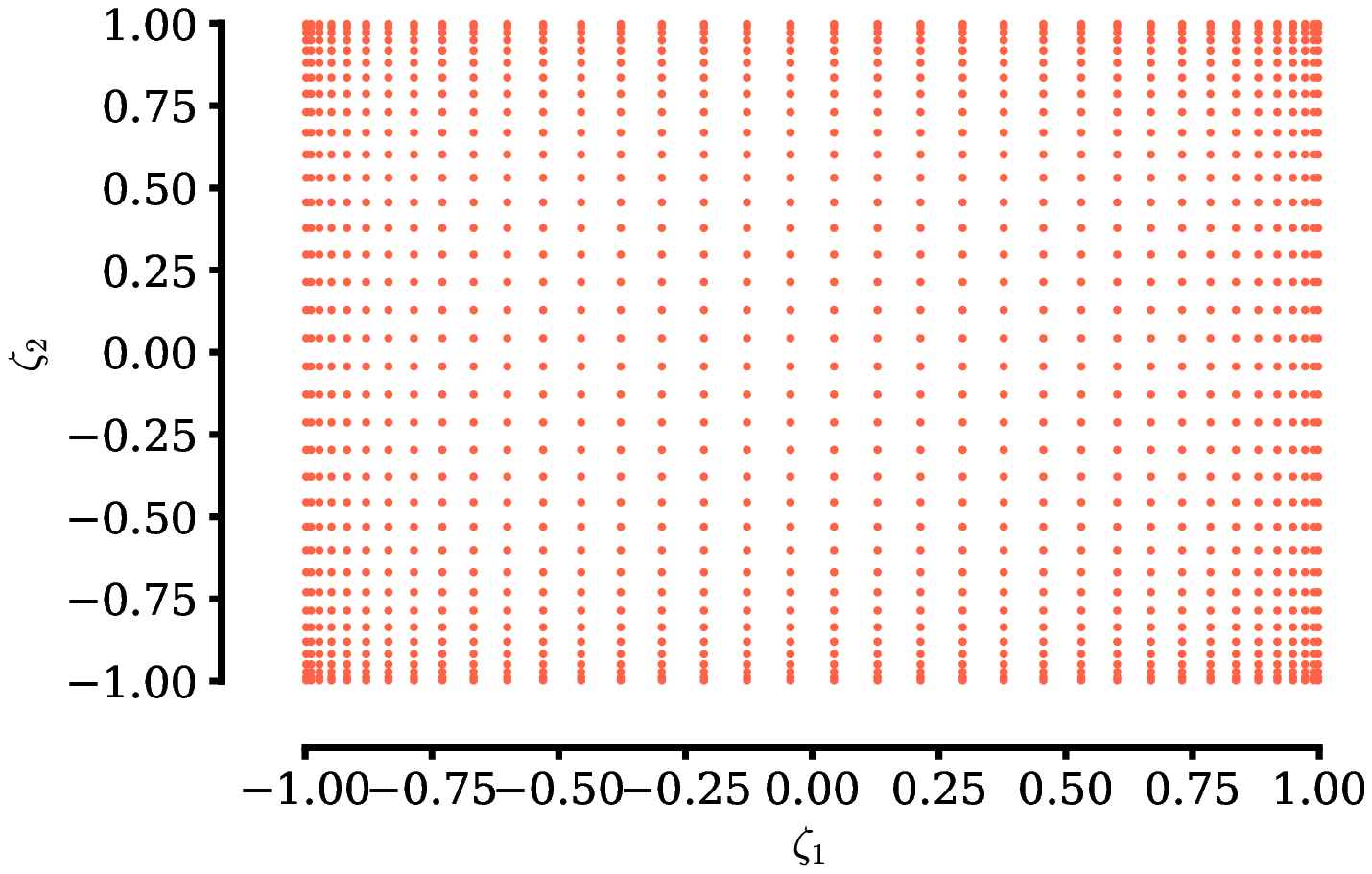}}
\subfigure[]{\includegraphics[]{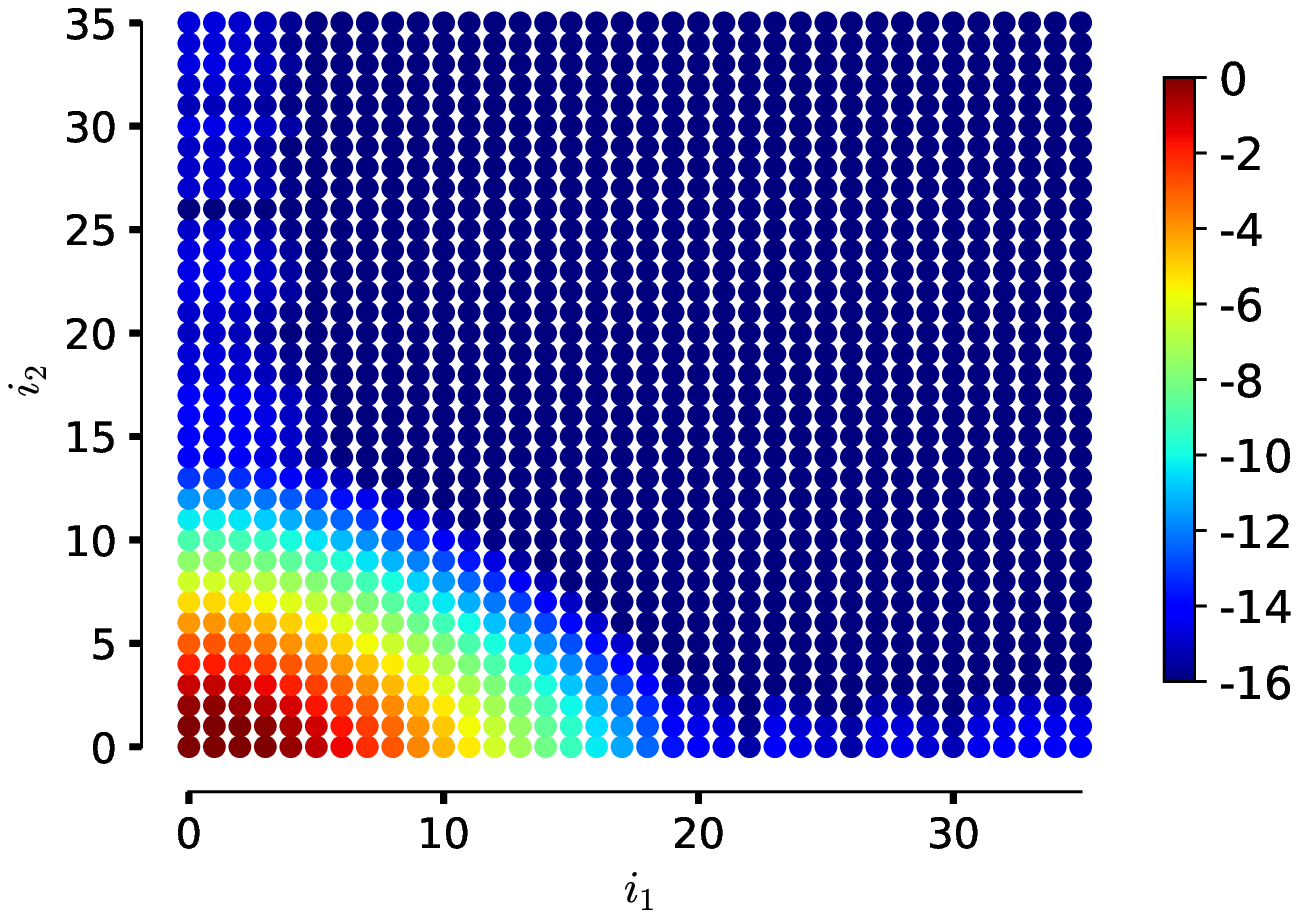}}
\subfigure[]{\includegraphics[]{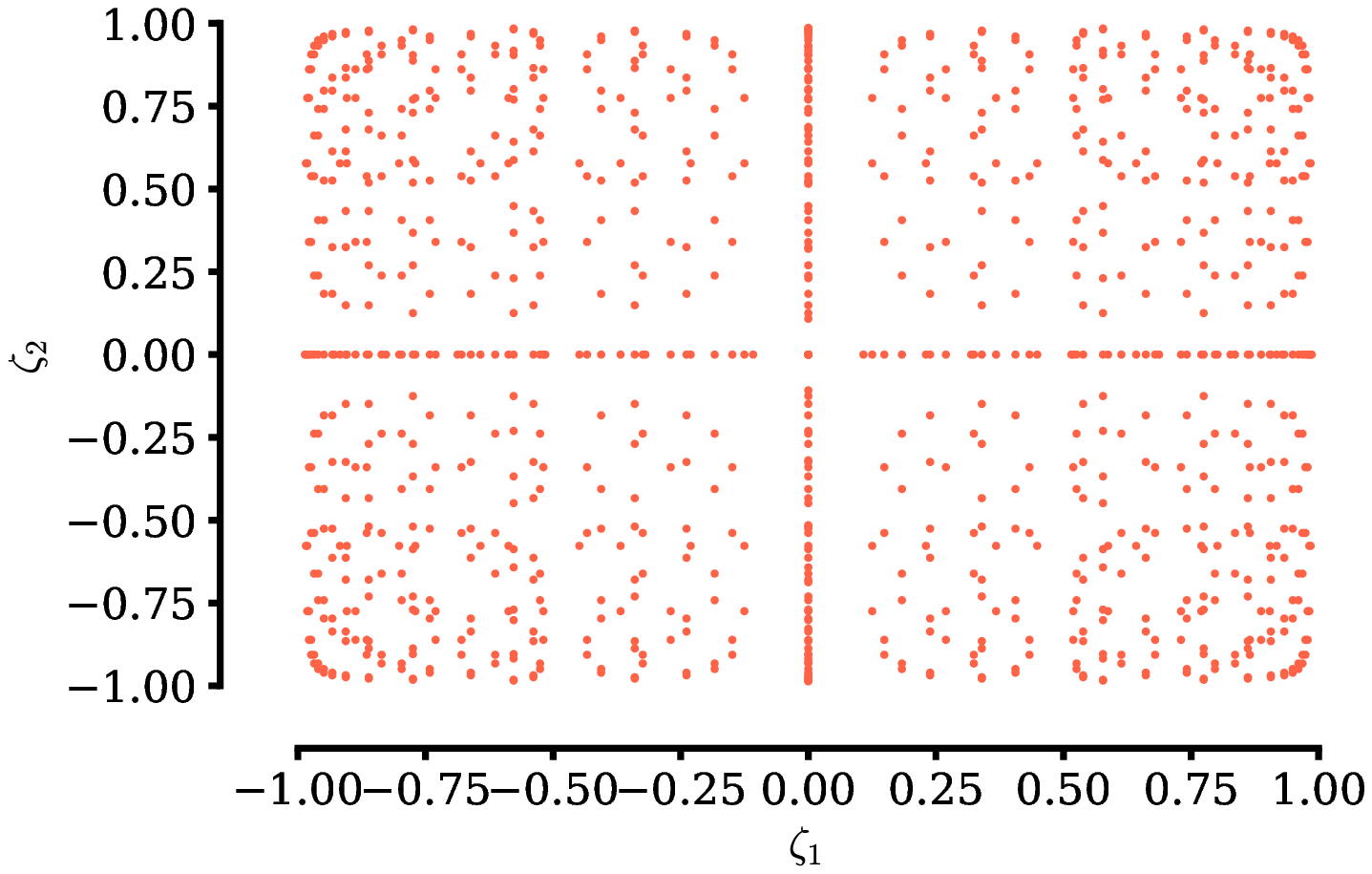}}
\subfigure[]{\includegraphics[]{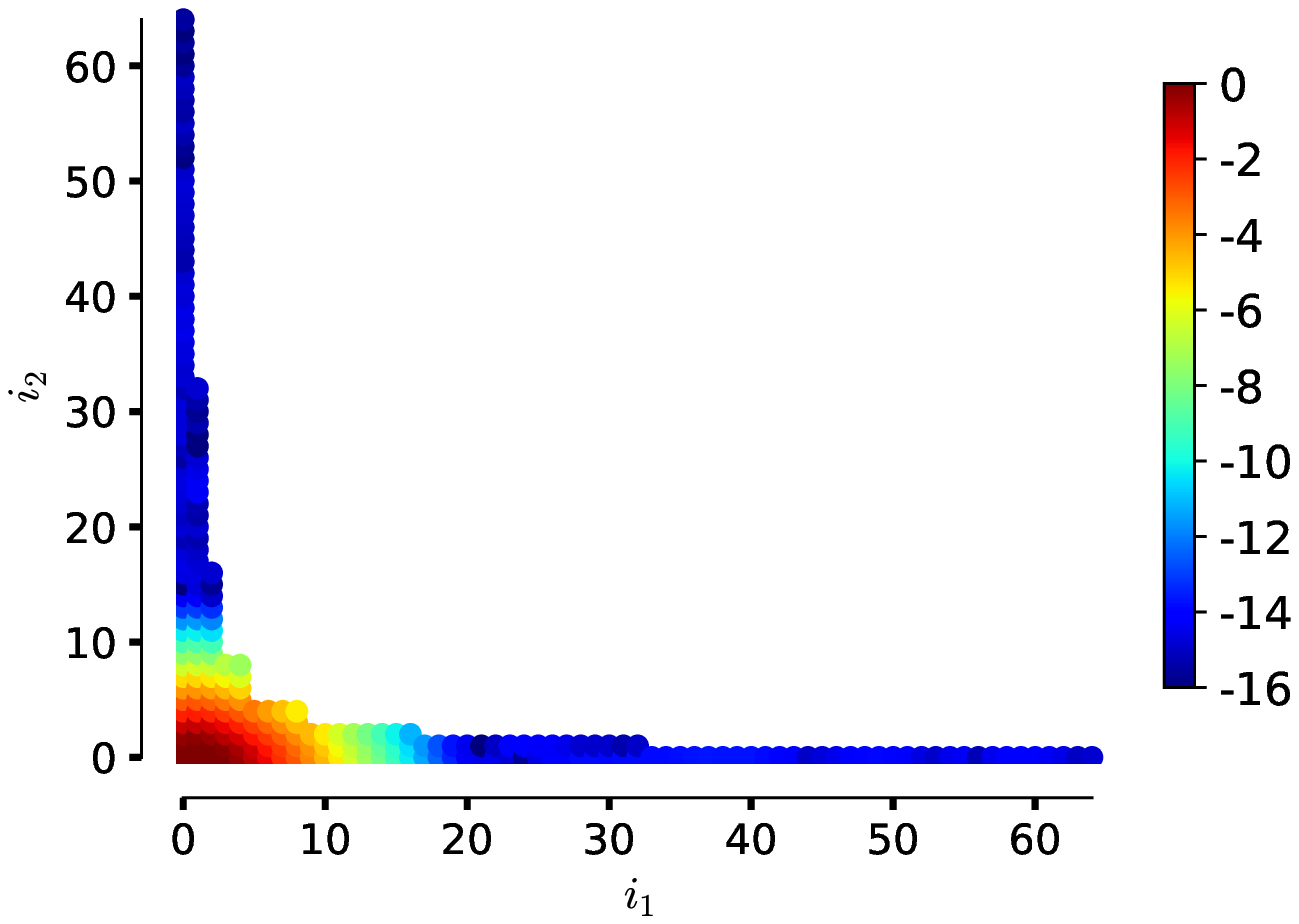}}
\subfigure[]{\includegraphics[]{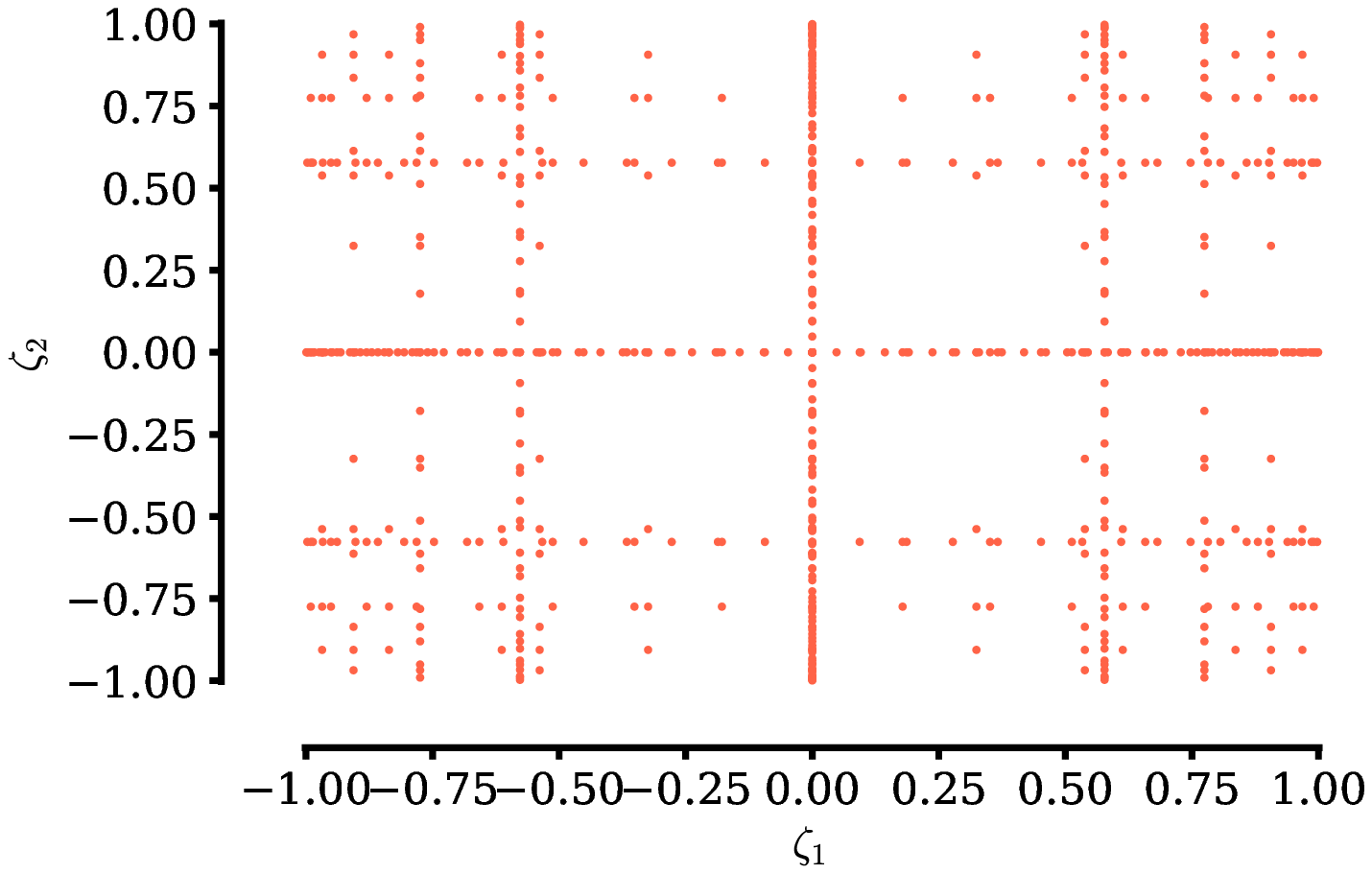}}
\subfigure[]{\includegraphics[]{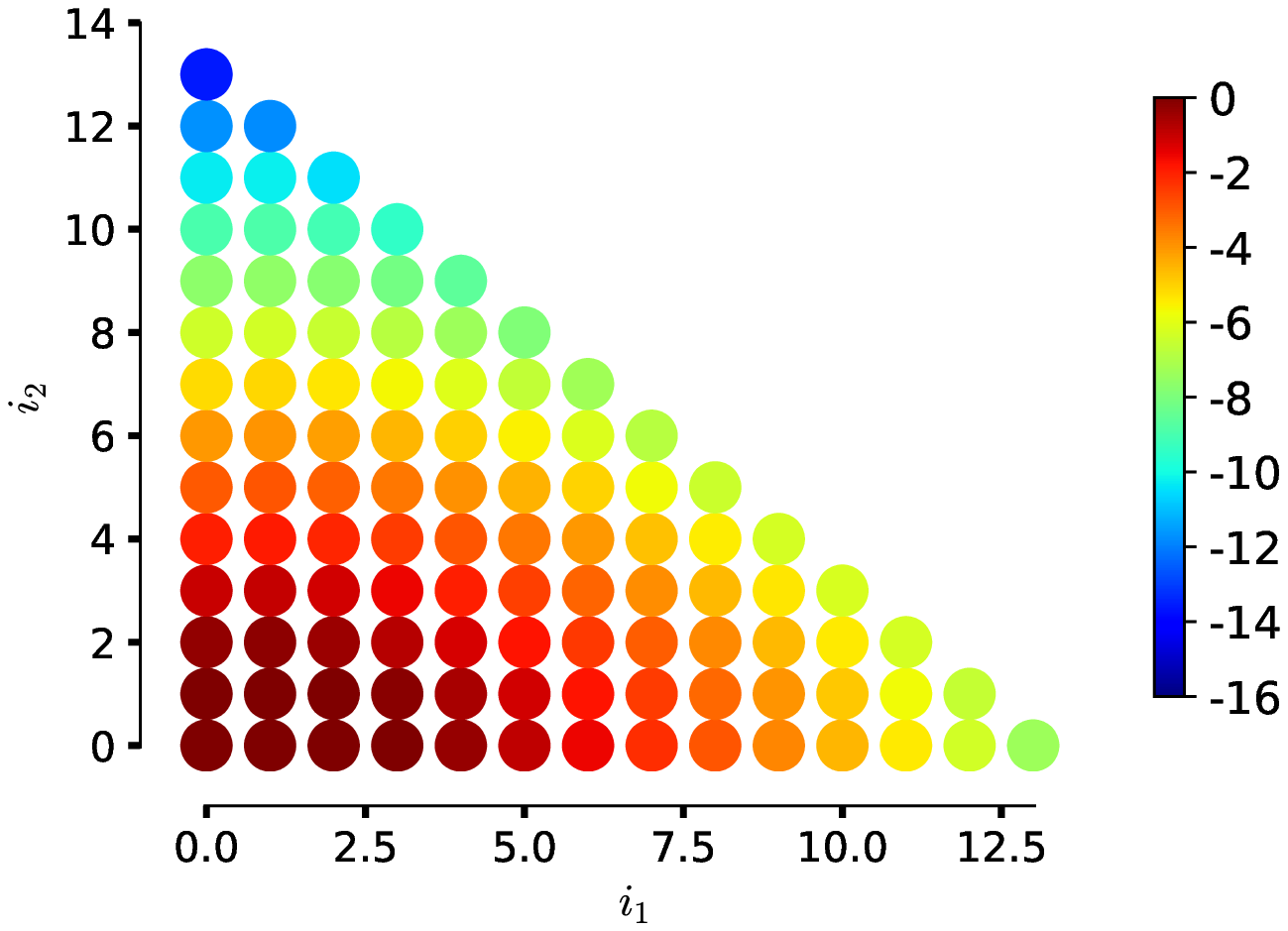}}
\end{subfigmatrix}
\caption{Parameter space points and pseudospectral coefficients (shown on a log base 10 scale) obtained via a (a, b) tensor grid with a maximum univariate order of 35; (c, d) sparse grid with a linear growth rule and $l=12$; (e, f) sparse grid with an exponential growth rule and $l=6$. The number of unique quadrature points from the three experiments were 1296, 1015 and 667 respectively.}
\label{sparsegrids}
\end{figure}

\begin{svgraybox}
\textbf{Approximation via sparse grids:}
To motivate the use of sparse grids, and to realize its limitations, consider the bi-variate function 
\begin{equation}
f\left( \bm{ \zeta}  \right)=exp\left(3 \zeta_{1}+\zeta_{2}\right)\; \; \; \text{where} \;  \; \; \bm{\zeta}\in\left[-1,1\right]^{2}
\end{equation}
and $\rho(\bm{\zeta})$ is the uniform distribution over the domain. Figure~\ref{sparsegrids}(a) plots the points in the parameter space for an isotropic tensor product grid using tensor products of order 35 Gauss-Legendre quadrature points, while (b) plots estimates of the pseudospectral coefficients using a tensor product Legendre polynomial basis in $\bm{\psi}_{\bm{p}}$. A total of 1296 function evaluations were required for these results. 

Upon inspecting (b) it is apparent that coefficients with a total order of 12 and greater can be set to zero with near negligible change in the function approximation. In Figure~\ref{sparsegrids}(c, d) and (e, f) we present two solutions that achieve this. Figure~\ref{sparsegrids}(c-d) plot the results for a sparse grid with a linear growth rule (essentially emulating a total order index set) with $l=12$, requring 1015 function evaluations, while in Figure~\ref{sparsegrids}(e-f) the results show an exponential growth rule with $l=5$, requiring 667 function evaluations. Clearly, using sparse grid quadrature rules can reduce the number of model evaluations compare to tensor product grids. 
\end{svgraybox}

While sparse grids and their adaptive variants \cite{conrad2013adaptive, pfluger2010spatially} are more computationally tractable than tensor product grids, they are still restricted to very specific index sets, even with linear, exponential and slow exponential \cite{burkardt2014slow} growth rules. Furthermore, when simulation-based function evaluations at the quadrature points fail (or are corrupted), one may have to resort to interpolation heuristics. In comparison, least squares methods offer far more flexibility---both in terms of a choice of the basis and in negotiating failed simulations.

\subsection{Scope and outline of paper}
Our specific goal in this paper is to use ideas from polynomial least squares to generate quadrature rules. Without loss of generality, these quadrature rules will be used to estimate the pseudospectral coefficients $\vx$ by solving
\begin{align}
\underset{\vx \in \mathbb{R}^{n}}{\text{minimize}}  \; \; \;  \left\Vert \mA \vx- \vb\right\Vert _{2},
 \label{leastsquares}
\end{align}
where we define the elements of $\mA \in \mathbb{R}^{m \times n}$ as per \eqref{def_of_A}. For completeness, we restate this definition but now in the multivariate setting, i.e., 
\begin{equation}
\mA(i, \bm{j} ) = \bm{\psi}_{\bm{j}} \left(\bm{\zeta} \right) \sqrt{\omega_{i}},  \; \; \; \; \; \bm{j} \in \mathcal{I}, \; \; \; \text{card}\left(\mathcal{I} \right) = n, \; \; \; i=1, \ldots, m.
\end{equation}
Typically $\mathcal{I}$ will be either a total order or a hyperbolic basis. The entries of the vector $\vb \in \mathbb{R}^{m}$ comprise of weighted model evaluations at the $m$ quadrature points; individual entries are given by $\vb(i) = \sqrt{\omega}_{i} f\left( \bm{\zeta}_{i} \right)$. Assuming $\mA$ and $\vb$ are known, solving \eqref{leastsquares} is trivial---the standard QR factorization can be used. However, the true challenge lies in selecting multivariate quadrature points and weights $\left\{ \left(\bm{\zeta}_{i},\omega_{i}\right)\right\} _{i=1}^{m}$, such that:
\begin{itemize}
\item For a given choice of $n$, the number of quadrature points $m$ can be minimized and yield a rule that has a high degree of exactness;
\item The least squares approximation in \eqref{leastsquares} is accurate;
\item The least squares solution is stable with respect to perturbations in $\vb$.
\end{itemize}
We will explore different strategies for generating $\mA$ and techniques for subsampling it---even introducing a new approach in secti.on~\ref{sec:convex}. Broadly speaking, strategies for computing multivariate quadrature points and weights via least squares involve two key decisions:
\begin{enumerate}
    \item \emph{Selecting a sampling strategy}: A suitable discretization of the domain from which quadrature points need to be computed; 
    \item \emph{Formulating an optimization problem}: The strategy for subselecting points from this sampling (if required) via the optimization of a suitable objective.
\end{enumerate}
Our goal in this manuscript is to describe the various techniques for generating these quadrature points and to make precise statements (where permissible) on their computational complexity. We substantiate our detailed review of literature with examples using the open-source code Effective Quadratures \cite{seshadri2017effective} \footnote{The codes to replicate the figures in this paper can be found at the website: \url{www.effective-quadratures.org/publications}}.

\section{Selecting a sampling strategy}
In this section we present methods for discretizing the domain based on the support of the parameters and their distributions. The survey builds on a recent review of sampling techniques for polynomial least squares  \cite{hadigol2017least} and from more recent texts including Narayan \cite{narayan2017computation, narayan2017christoffel} and Jakeman and Narayan \cite{jakeman2017generation}. 

Intuitively, the simplest sampling strategy for polynomial least squares is to generate random Monte-Carlo type samples based on the joint density $\rho(\bm{\zeta})$. Migliorati et al. \cite{migliorati2014analysis} provide a practical sampling heuristic for obtaining well conditioned matrices $\mA \in \mathbb{R}^{m \times n}$; they suggest sampling with $m=n^2$ for total degree and hyperbolic cross spaces. Furthermore, they prove that as the number of samples $m \rightarrow \infty$, the condition number $\kappa \left( \mA^T \mA \right) \rightarrow 1$ (see page 441 in \cite{migliorati2014analysis}). While these random samples can be pruned down using the optimization strategies presented in~\ref{sec:opt}, there are other sampling alternatives.

\subsection{Christoffel samples}
\label{sec:christoffel}
The Christoffel sampling recipe of Narayan and co-authors \cite{narayan2017christoffel} significantly outperforms the Monte Carlo approach in most problems. To understand their strategy, it will be useful to define the diagonal of the orthogonal projection kernel 
\begin{equation}
K_n \left( \bm{\zeta} \right) \coloneqq \sum_{ j = 1 }^{n} \bm{\psi}_{\bm{j} }^{2}\left( \bm{\zeta}  \right),
\end{equation}
where as before the subscript $n$ denotes $\text{card} \left(\mathcal{I} \right)$. We associate the above term with a constant 
\begin{equation}
\left\Vert K\right\Vert _{\infty}\triangleq\underset{\bm{\zeta} \in\mathbb{R}^{D}}{sup}\left(K_{n}\left( \bm{\zeta} \right)\right),
\end{equation}
which has a lower bound of $n$. Cohen et al. \cite{cohen2013stability} prove that if the number of samples $m$ satisfies the bound 
\begin{equation}
\frac{m}{nlog\left(m \right)}\apprge\chi\frac{\left\Vert K\right\Vert _{\infty}}{n},
\label{sampling_bound}
\end{equation}
for some constant $\chi$, then the subsequent least squares problem is both stable and accurate with high probability. The key ingredient to this bound, the quantity $\left\Vert K\right\Vert _{\infty}$, tends to be rather small for Chebyshev polynomials---scaling linearly with $n$ in the univariate case---but will be large for other polynomials, such as Legendre polynomials---scaling quadratically with $n$. Multivariate examples are given in \cite{chkifa2015discrete} and lead to computationally demanding restrictions on the number of samples required to satisfy \eqref{sampling_bound} \cite{narayan2017christoffel, cohen2017optimal}. Working with this bound, Narayan et al. provide a sampling strategy that leverages the fact that the asymptotic characteristics for total order polynomials can be determined. Thus, if the domain is bounded with a continuous $\rho$, then limit of $m/K_n$---known as the Christoffel function---is given by
\begin{equation}
\underset{n\rightarrow\infty}{ \text{lim} } \; \frac{m}{K_{n}\left( \bm{\zeta} \right)}=\frac{ \bm{\rho} \left(\bm{\zeta}\right)}{\nu\left(\bm{\zeta}\right)}
\end{equation}
where the above statements holds weakly, and where $\nu (\bm{\zeta})$ is the \emph{weighted pluripotential equlibrium measure} (see section 2.3 in \cite{narayan2017christoffel} for the definition and significance). This implies that the functions $\bm{\phi} = \bm{\psi} \sqrt{m / K_n}$ form an orthogonal basis with respect to a modified weight function $\bm{\rho} K_n / m$. This in turn yields a modified reproducing diagonal kernel
\begin{equation}
\hat{K}_{n}=\sum_{j=1}^{m}  \bm{\phi}^{2}=\frac{n}{K_{n}}K_{n}=n,
\end{equation}
which attains the optimal value of $n$. The essence of the sampling strategy is to sample from $\nu$ (performing a Monte Carlo on the basis $\bm{\phi}$), which should in theory reduce the sample count as dictated by \eqref{sampling_bound}. 

The challenge however is devleoping computational algorithms for generating samples according to $\nu(\bm{\zeta})$, since precise analytical forms are only known for a few domains. For instance, when $D = [-1,1]^d$, the measure  $\nu(\bm{\zeta})$ is the Chebyshev (arcsine) distribution. Formulations for other distributions such as the Gaussian and exponential distributions are in general more complex and can be found in section 6 of \cite{narayan2017christoffel}.

\begin{figure}
\begin{subfigmatrix}{2}
\subfigure[]{\includegraphics[]{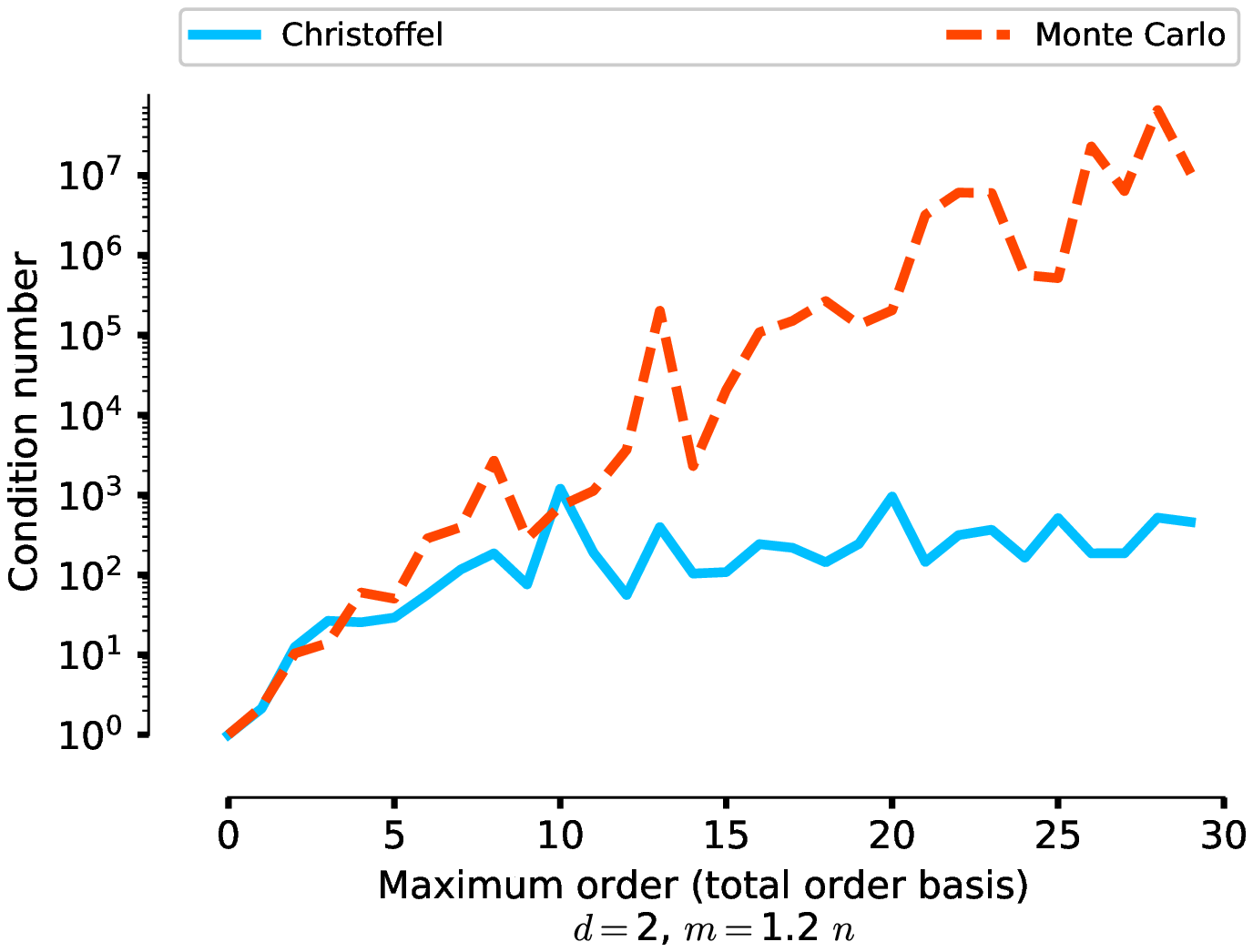}}
\subfigure[]{\includegraphics[]{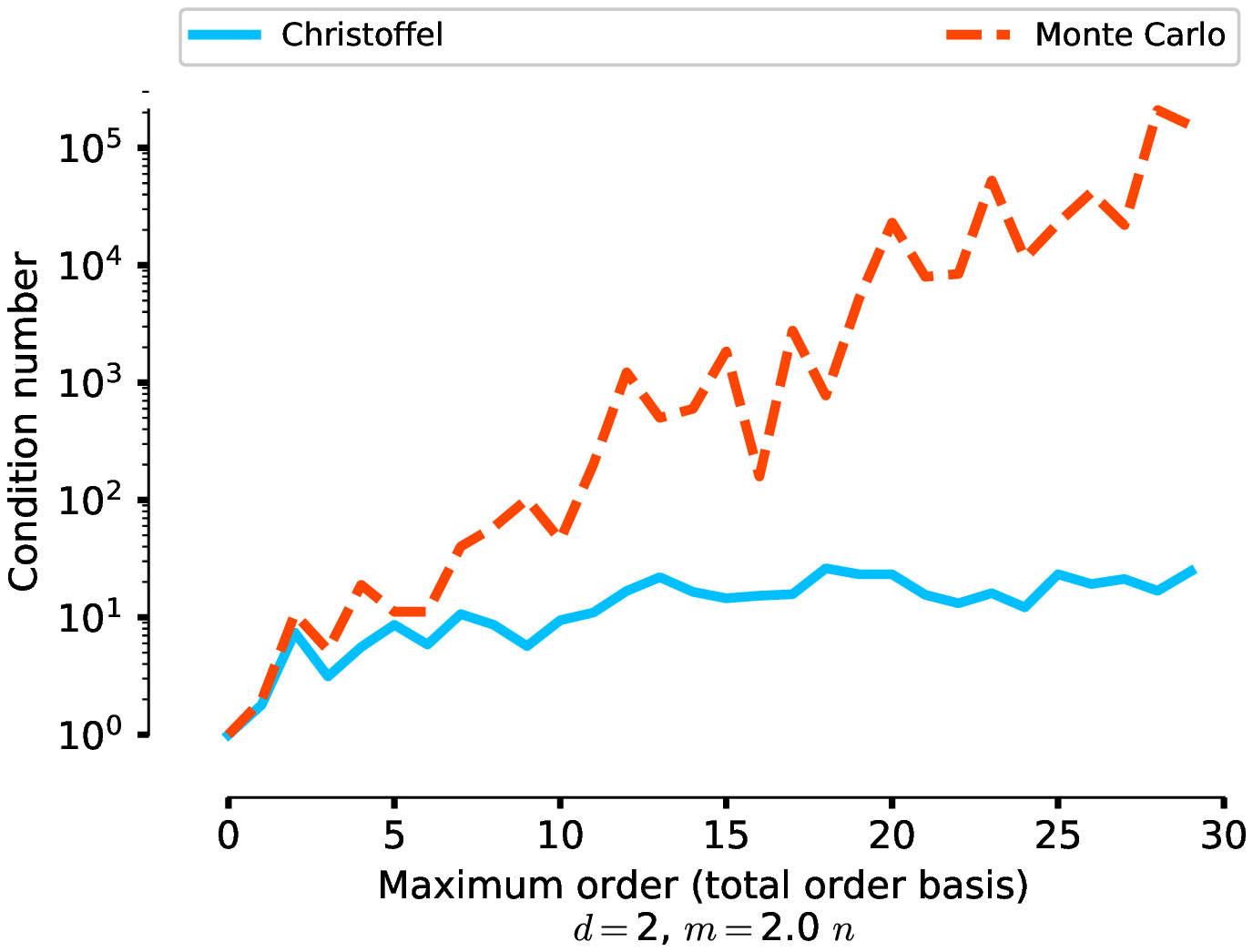}}
\subfigure[]{\includegraphics[]{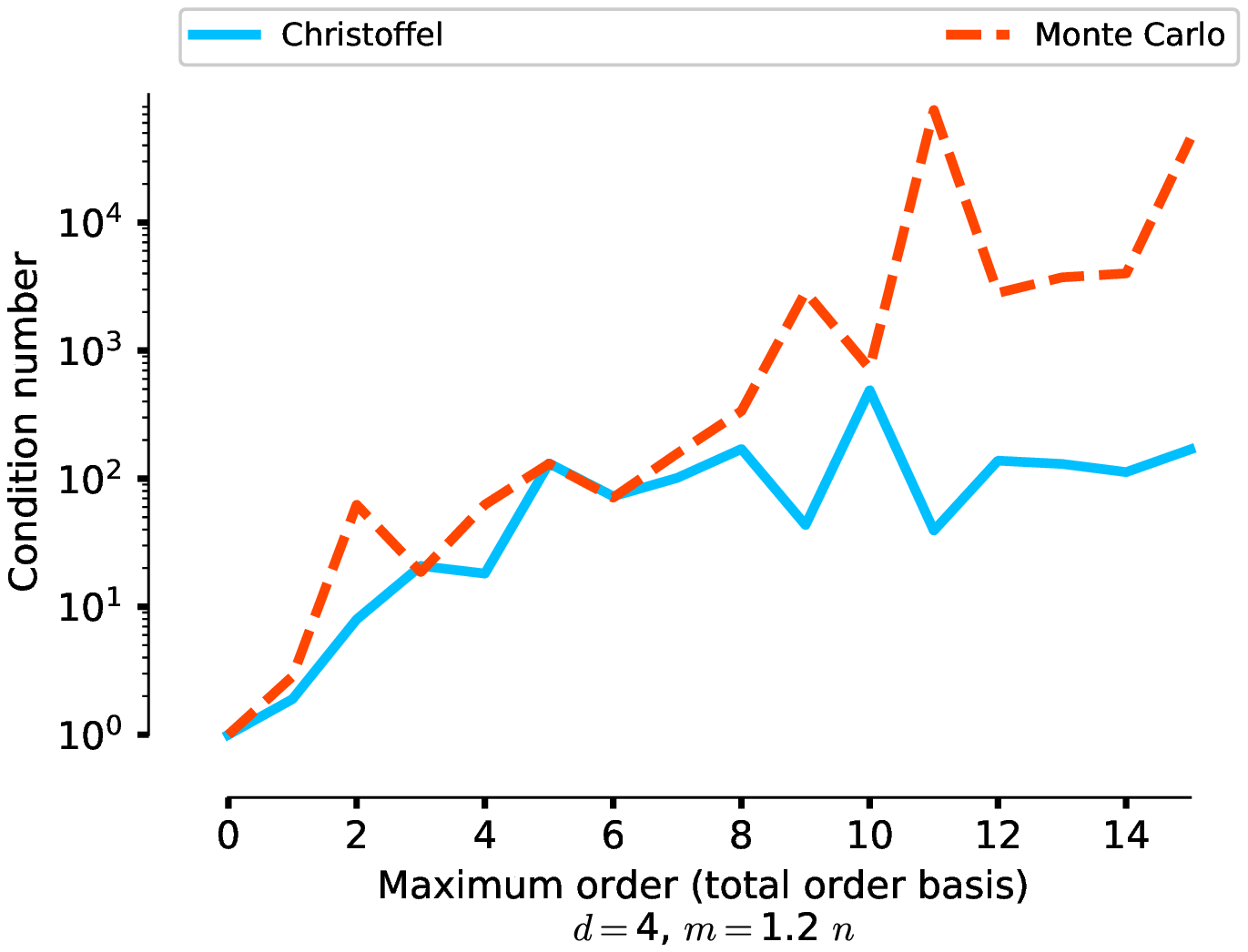}}
\subfigure[]{\includegraphics[]{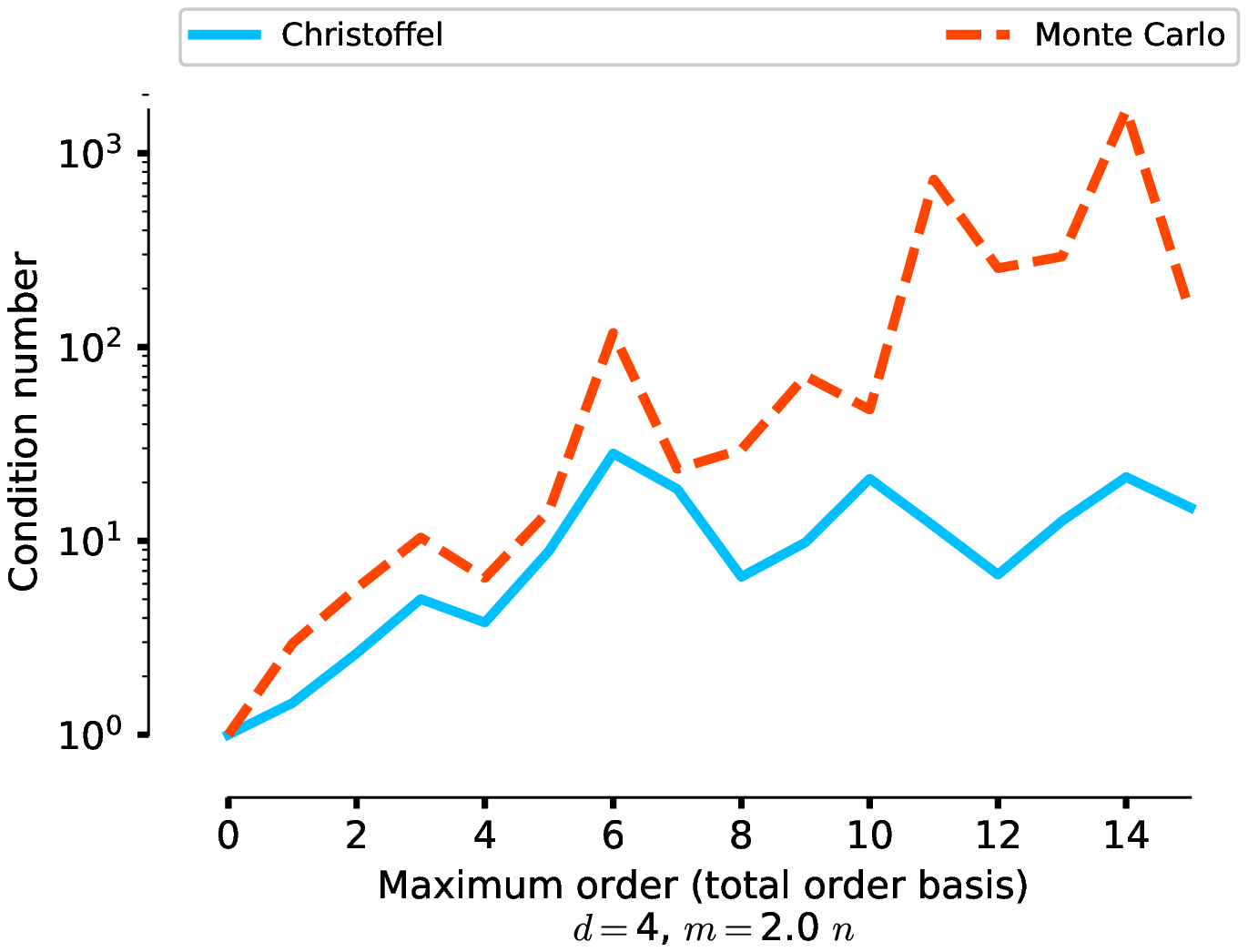}}
\end{subfigmatrix}
\caption{A comparison of condition numbers for $\mA$ when constructing the matrix using multivariate Legendre polynomials evaluated at either Monte Carlo or Christoffel samples. Experiments for (a) $d=2, m=1.2n$; (b) $d=2, m=2n$; (c) $d=4, m=1.2n$; (d) $d=2, m=2n$. Plotted are the mean values of 10 random trials.}
\label{cs_demo}
\end{figure}

\begin{figure}
\centering
\includegraphics[scale=0.4]{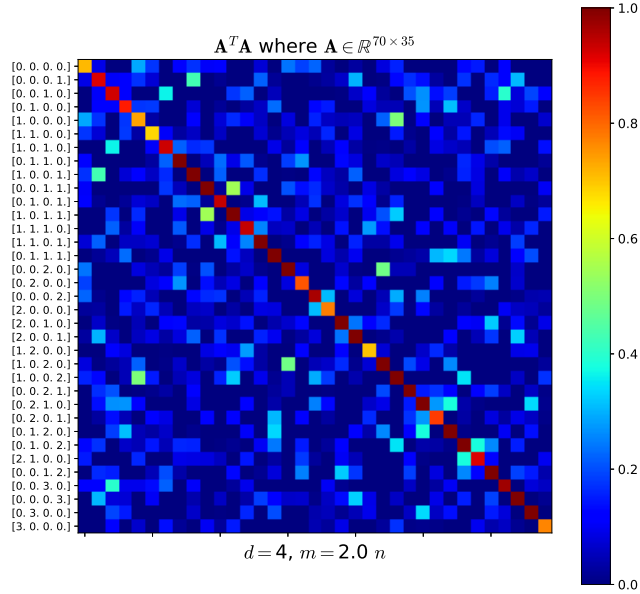}
\caption{The Gram matrix $\mG = \mA^{T} \mA$ for one of the trials for the case where $d=4$ and $m=2n$ showing its deviation from the identity.}
\label{cs_demo_2}
\end{figure}

\begin{svgraybox}
\textbf{Comparing condition numbers:}
But how much lower are the condition numbers when we compare standard Monte Carlo with Christoffel sampling? Figure~\ref{cs_demo} plots the mean condition numbers (averaged over 10 trials) for Legendre polynomials in $d=2$ in (a) and (b), and $d=4$ in (c) and (d). Two different oversampling factors are also applied; (a) and (c) have an oversampling factor of $1.2$ while (b) and (d) have an oversampling factor of $2$--i.e., $m=2n$. For the Christoffel results, the samples are generated from the Chebyshev distribution. 

It is apparent from these figures that the Christoffel sampling strategy does produce more well conditioned matrices on average. We make two additional remarks here. The first concerns the choice of the weights. In both sampling strategies the \emph{quadrature weights} $\omega_{i}$ are computed via
\begin{equation}
\omega_{i}=\frac{\tilde{\omega}_{i}}{\sum_{k=1}^{m}\tilde{\omega}_{k}} ,\; \; \; \text{where} \; \; \;  \tilde{\omega}_{i}=\frac{n}{m}\sum_{\bm{j} \in\mathcal{I}} \bm{\psi}_{\bm{j}}^{2}\left( \bm{\zeta}_{i}\right)
\end{equation}
ensuring that the weights sum up to unity.  The second concerns the Gram matrix $\mA^{T} \mA$, shown in Figure~\ref{cs_demo_2} for the Christoffel case with a maximum order of 3, with $d=4$ and an oversampling factor of $2$. Although the condition number of this matrix is very low (4.826), one can clearly observe sufficient internal aliasing errors that would likely effect subsequent numerical computation of moments. 
\end{svgraybox}

\subsection{Subsampling tensor grids}
\begin{figure}
\begin{subfigmatrix}{2}
\subfigure[]{\includegraphics[]{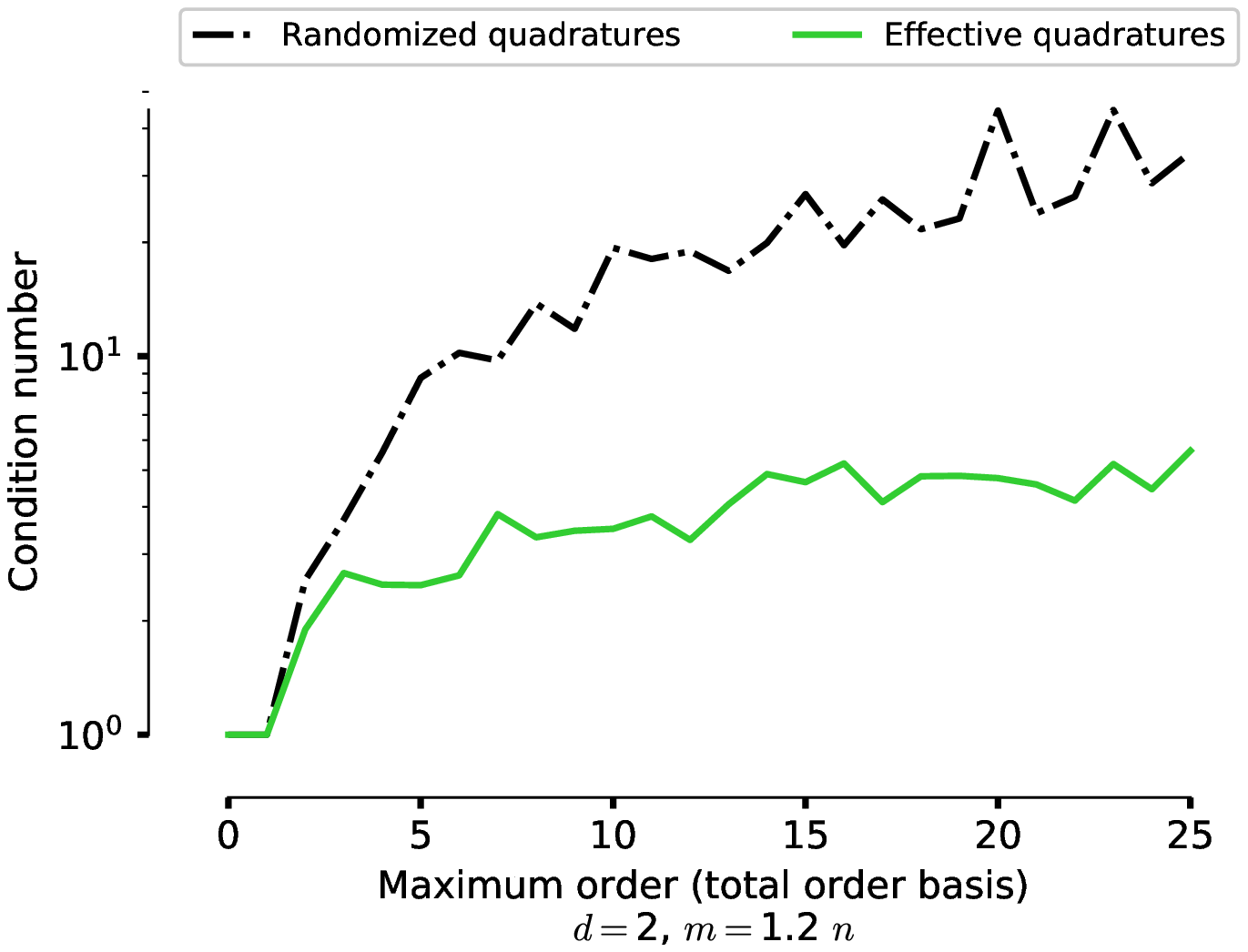}}
\subfigure[]{\includegraphics[]{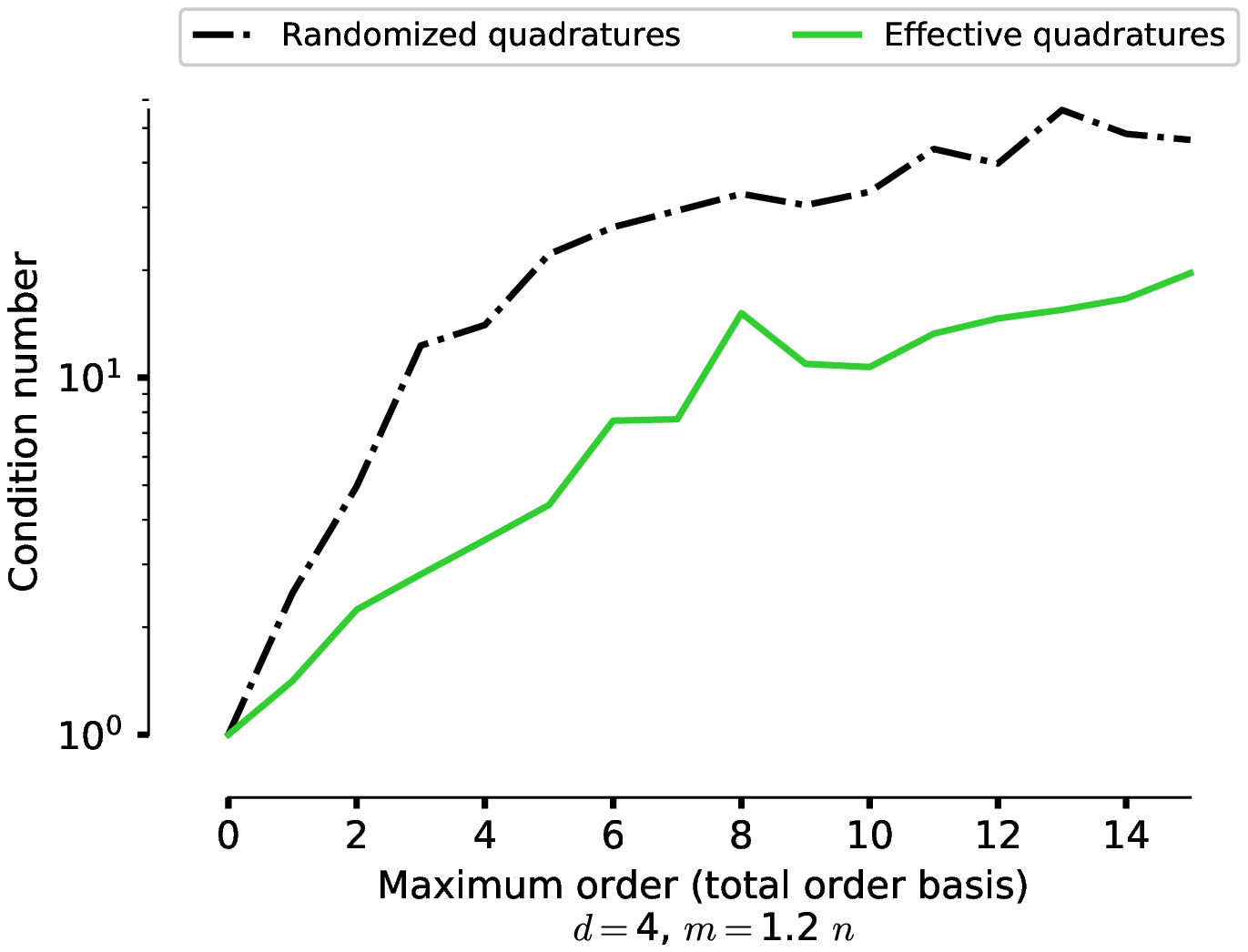}}
\end{subfigmatrix}
\caption{A comparison of condition numbers for $\mA$ when constructing the matrix using multivariate Legendre polynomials evaluated at either randomized or effectively subsampled quadrature points. Experiments for (a) $d=2, m=1.2n$; (b) $d=4, m=1.2n$. For the randomized technique, we plot the mean values of 10 random trials.}
\label{eqrq_demo}
\end{figure}

The idea of constructing $\mA$ using samples from a tensor grid and then subsampling its rows has been explored from both compressed sensing\footnote{By solving the basis pursuit (and de-noising) problem.} \cite{tang2014subsampled} and least squares \cite{zhou2015weighted, seshadri2017effectively} perspectives. In Zhou et al. \cite{zhou2015weighted} the authors randomly subsample the rows and demonstrate stability in the least squares problem with $m$ scaling linearly with $n$. In Seshadri et al., \cite{seshadri2017effectively}, a greedy optimization strategy is used to determine which subsamples to select; they report small condition numbers even when $n=m$ for problems where $d \leq 7$ and a total order of degree $ \leq 15$. Details of their optimization strategy are presented in section~\ref{sec:opt}. One drawback of their technique is the requirement that the full $\mA$ matrix must be stored and passed onto the optimizer---a requirement that can be circumvented in the \emph{randomized} approach, albeit at the cost of $m$ being greater than $n$ (typically). A representative comparison of the condition numbers is shown in Figure \ref{eqrq_demo}. In these results, the isotropic tensor grid from which the subsamples are computed corresponds to the highest total order of the polynomial basis (shown on the horizontal axis).

\subsection{Coherence and induced sampling} 
Building on some of the ideas discussed previously in \ref{sec:christoffel}, one can detail an \emph{importance sampling} strategy that yields stable least squares estimates for $m \apprge n$, up to log factors. The essence of the idea is to define a new distribution $\mu$ such that
\begin{equation}
\mu=\frac{K_{n}}{n} \bm{\rho}=\frac{1}{n}\left(\sum_{j=1}^{n} \bm{\psi}_{j}^{2}\left( \bm{\zeta} \right)\right) \bm{\rho}.
\end{equation}
It should be noted that while $\bm{\rho}$ is a product measure and therefore easy to sample from, $\mu$ is not and thus requires either techniques based on Markov chain Monte Carlo (MCMC) or conditional sampling \cite{cohen2017optimal}. More specifically, $\mu$ depends on $m$, so if one enriches the space with more samples or change the basis, then $\mu$ will change. In \cite{hampton2015coherence}, the authors devise an MCMC strategy which seeks to find $\mu$, thereby explicitly minimizing their \emph{coherence} parameter (a weighted form of $\left\Vert K\right\Vert _{\infty}$ ). In \cite{narayan2017computation}, formulations for computing $\mu$ via induced distributions is detailed. In our numerical experiments investigating the condition numbers of matrices obtained via such \emph{induced samples}---in a similar vein to the Figures~\ref{cs_demo}---the condition numbers were found to be comparable to those from Christoffel sampling.


\section{Optimization strategies}
\label{sec:opt}
In the case of Christoffel, Monte Carlo and even subsampling based techniques, a reduction in $m$ can be achieved, facilitating near quadrature like degree of exactness properties. In this section we discuss various optimization techniques for achieving this.

\subsection{Greedy linear algebra approaches}
It will be convenient to interpret our objective as identifying a suitable submatrix of $\mA \in \mathbb{R}^{m \times n}$ by only selecting $k <<m$ rows. Ideally, we would like $k$ to be as close to $n$ as possible. Formally, we write this submatrix as
\begin{align}
\begin{split}
\hat{\mA}_{\vz} &=  \frac{1}{\tau} \sum_{i=1}^{m} z_{i} \va_{i}^T \ve \\
\text{where}  \; \; \; & \mathbf{1}^{T} \vz =k, \\
 & \tau = \sum_{i=1}^{m}z_{i}w_{i}, \\
 & z_{i}\in\left\{ 0,1\right\}
 \end{split}
 \label{cvx}
\end{align}
and where $\va_{i}^{T}$ represents the rows of $\mA$, and $\ve \in \mathbb{R}^{k}$ and $\mathbf{1} \in \mathbb{R}^{m}$ are vector of ones. The notation above indicates that $\vz \in \mathbb{R}^{m}$ is a boolean vector. The normalization factor $\tau$  is introduced to ensure the weights associated with the \emph{subselected} quadrature rule sum up to one.

Table 1 highlights some of the various optimization objectives that may be pursued for finding $\vz$ such that the resulting $\hat{\mA}_{\vz}$ yields accurate and stable least squares solutions. 

\begin{table}
\begin{center}
\caption{Sample objectives for optimization}
\vspace{2 mm}
\begin{tabular}{l l }
\hline
\textbf{Case} & \textbf{Objective} \\
\hline
P1: & $\underset{\vz}{\text{minimize}} \; \; \sigma_{max} \left( \hat{\mA}_{\vz} \right)$ \\
P2: & $\underset{\vz}{\text{maximize}} \; \; \sigma_{min} \left( \hat{\mA}_{\vz} \right)$ \\
P3: & $\underset{\vz}{\text{minimize}} \; \; \kappa \left( \hat{\mA}_{\vz} \right)$ \\
P4:  & $\underset{\vz}{\text{maximize}} \; \; \text{vol} \left( \hat{\mA}_{\vz} \right) = \prod_{i=1}^{k}\sigma_{i} \left( \hat{\mA}_{\vz} \right)$ \\
P5: & $\underset{\vz}{\text{minimize}} \; \; \left\Vert \hat{\mA}_{\vz}^{T} \hat{\mA}_{\vz} \right\Vert _{2}$
\vspace{1 mm}\\
\hline
\end{tabular}
\end{center}
\label{table1}
\end{table}

We remark here that objectives P1 to P4 are proven NP-hard problems (see Theorem 4 in Civril and Magdon-Ismail \cite{ccivril2009selecting}); although it is readily apparent that all the objectives require evaluating $\binom{m}{k}$ possible choices---a computationally unwieldy task for large values of $m$ and $k$. Thus some regularization or relaxation is necessary for tractable optimization strategies.

We begin by focusing on P4. In the case where $k=n$, one can express the volume maximization objective as a determinant maximization problem. Points that in theory maximize this determinant, i.e., the \emph{Fekete points}, yield a Lebesgue constant that typically grows logarithmically---at most linearly---in $n$ \cite{bos2006bivariate}. So how do we find a determinant maximizing submatrix of $\mA$ by selecting $k$ of its rows? In Guo et al. \cite{guo2018weighted} the authors show that if there exists a point stencil that indeed maximizes the determinant, then a greedy optimization \emph{can} recover the global optimum (see Theorem 3.2 in Guo et al.). Furthermore, the authors prove that either maximizing the determinant or minimizing the condition number (P3) is likely to yield equivalent solutions. This explains their use of the pivoted QR factorization for finding such points; its factorization is given by
\begin{equation}
\mA^{T} \mP= \mQ\left(\begin{array}{cc}
\mR_{1} & \mR_{2}\end{array}\right),
\end{equation}
where $\mQ \in \mathbb{R}^{n \times n}$ is an orthogonal matrix and $\mR_{1} \in \mathbb{R}^{n \times n}$ is an upper triangular matrix that is invertible and $\mathbb{R}^{n \times (m-n)}$. The permutation matrix $\mP \in \mathbb{R}^{m \times m}$ is based on a pivoting strategy that seeks to either maximize $\sigma_{max}(\mR_{1})$ or minimize $\sigma_{min}(\mR_{2})$. As Bj\"{o}rck notes (see section 2.4.3 in \cite{bjorck2016numerical}) both these strategies are in some cases equivalent and thus greedily maximize the diagonal entries of $\mR_{1}$ and therefore serve as heuristics for achieving all the objectives in Table~\ref{table1}. Comprehensive analysis on pivoted QR factorizations can be found in \cite{chandrasekaran1994rank, chan1994low, dax2000modified}. Techniques for the closely related subset selection problem can also be adopted and build on ideas based on random matrix theory (see Deshpande and Rademacher \cite{deshpande2010efficient} and references therein). The monograph by Miller \cite{miller2002subset} which also provides a thorough survey of techniques, emphasizes the advantages of methods based on QR factorizations. 

\begin{figure}
\begin{subfigmatrix}{2}
\subfigure[]{\includegraphics[]{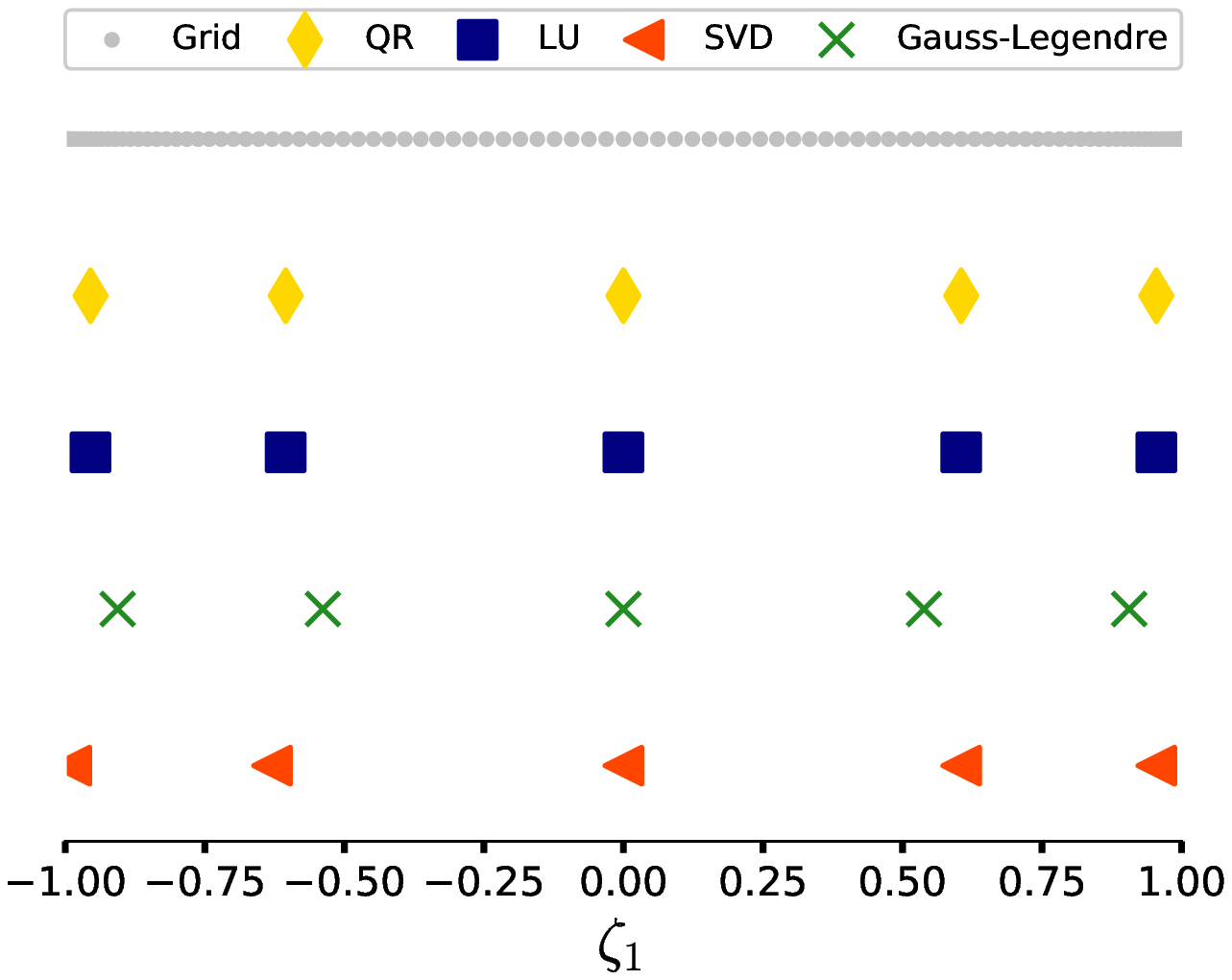}}
\subfigure[]{\includegraphics[]{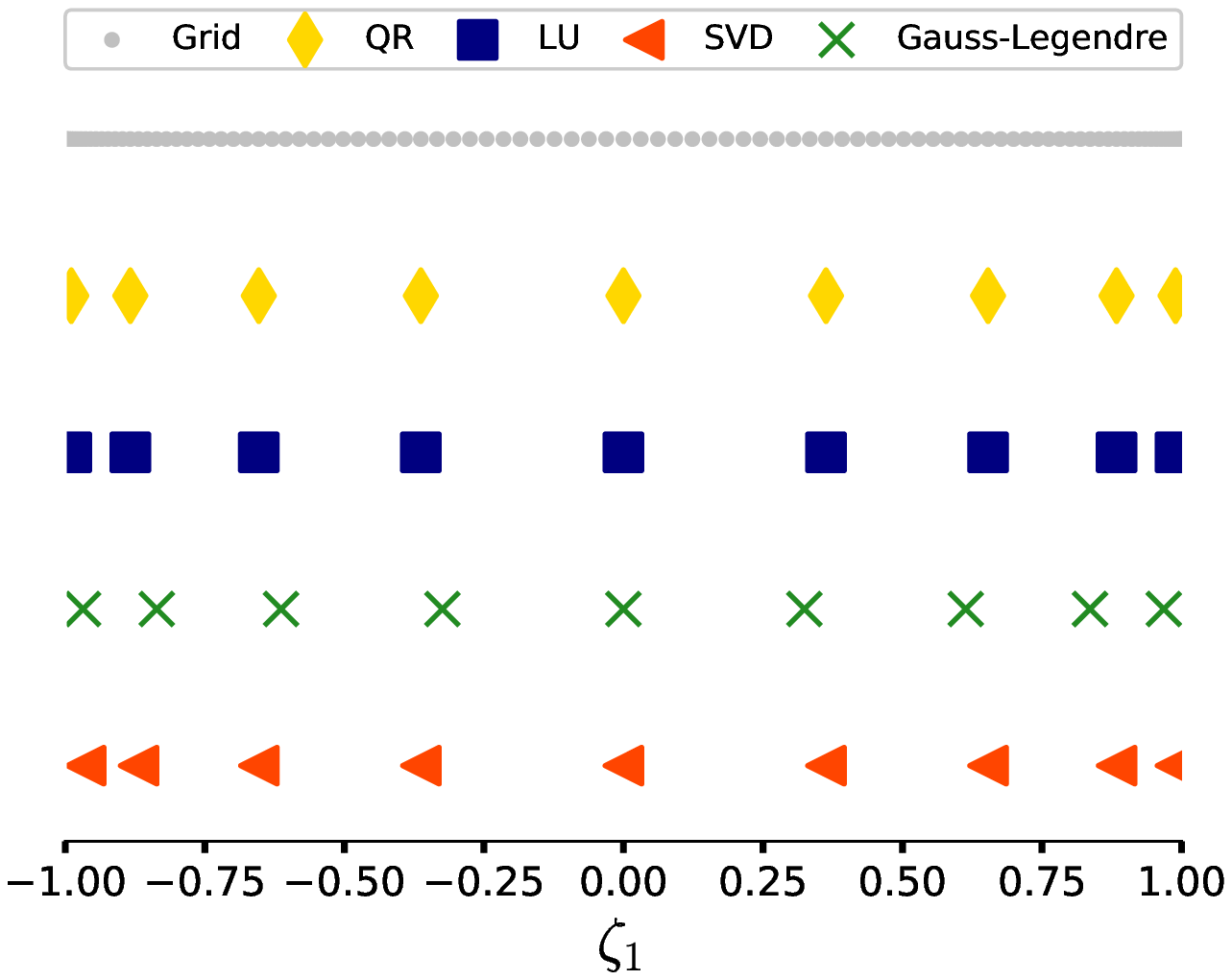}}
\end{subfigmatrix}
\caption{Approximation to the Gauss-Legendre quadrature points obtained by subsampling from a grid of order 100 using various optimization strategies with an (a) order 4; (b) order 8 Legendre polynomial basis.}
\label{gauss1d}
\end{figure}

\begin{figure}
\begin{subfigmatrix}{2}
\subfigure[]{\includegraphics[]{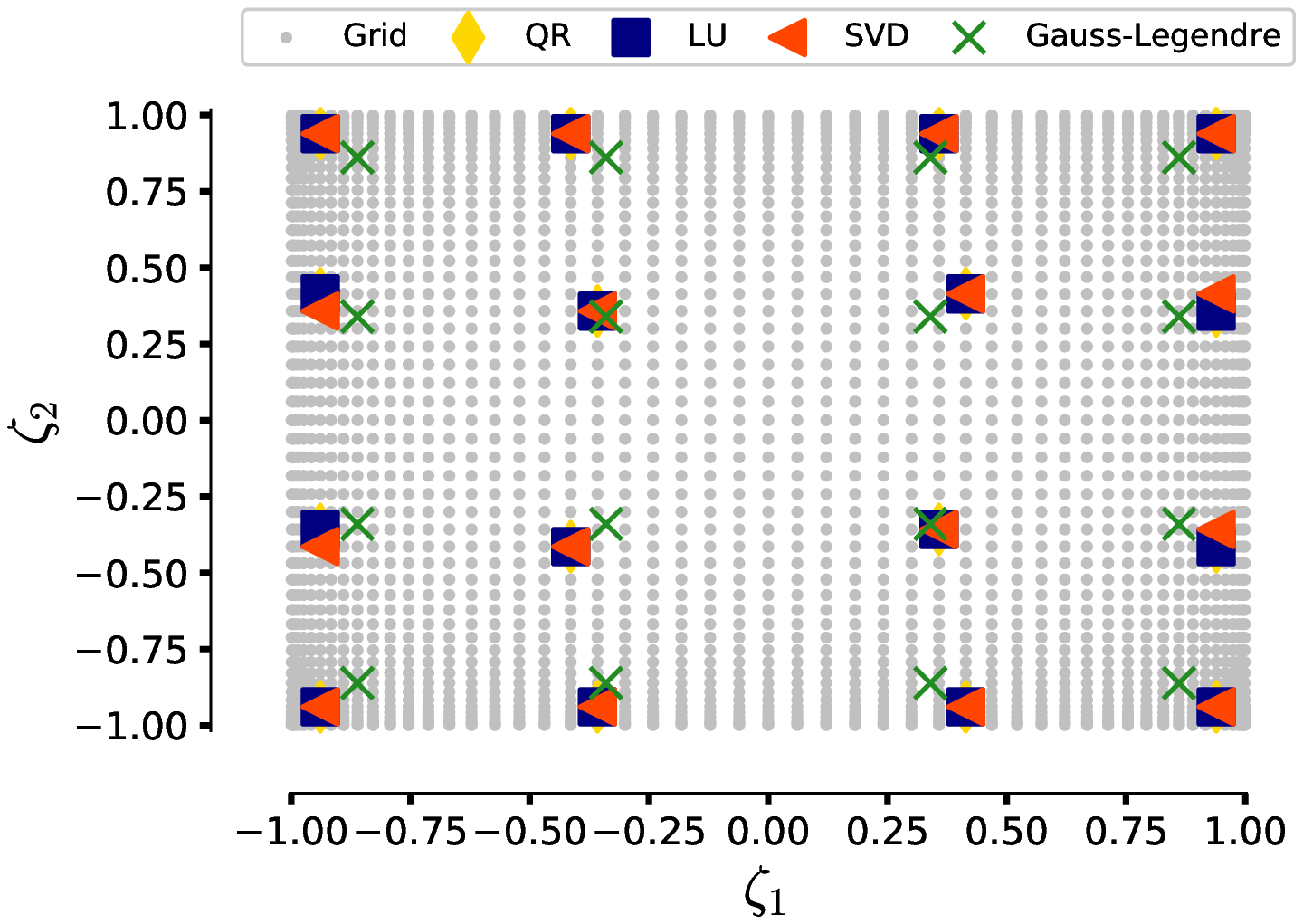}}
\subfigure[]{\includegraphics[]{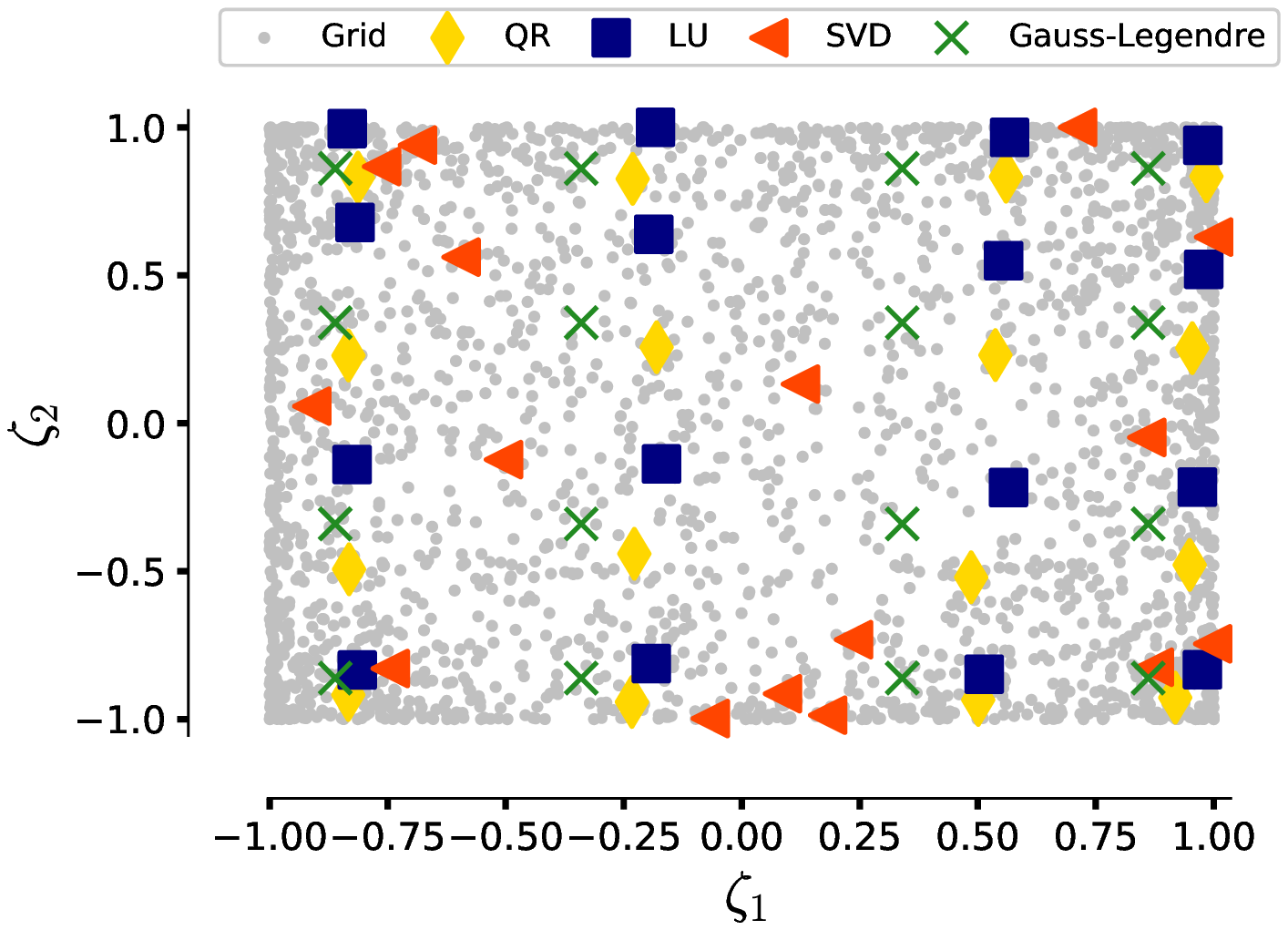}}
\end{subfigmatrix}
\caption{Approximation to the $[3, 3]$ Gauss-Legendre quadrature points obtained by subsampling from a grid of with 2601 samples using various optimization strategies: (a) Subsampling from a tensor grid of order $[50,50]$; (b) Subsampling from 2601 randomly distributed Chebyshev points.}
\label{gauss2d}
\end{figure}

\begin{svgraybox}
\textbf{Optimizing to find Gauss points:}
In the univariate case, it is known that Gauss-Legendre quadrature points are optimal with respect to a uniform measure. In an effort to gauge the efficacy of some of the optimization strategies discussed in this paper, we carry out a simple numerical experiment. Let $\mA \in \mathbb{R}^{101 \times K}$ be formed by evaluating up to order $K$ Legendre polynomials at the first 101 Gauss-Legendre quadrature points. Each of the aforementioned greedy strategies is tasked with finding a suitable submatrix $\hat{\mA}_{\vz} \in \mathbb{R}^{K \times K}$. The quadrature points embedded in $\vz$ for $K=4$ and $K=8$ are shown in Figure~\ref{gauss1d}(a) and (b) respectively. It is interesting to observe how the points obtained from LU with row pivoting, QR with column pivoting and SVD-based subset selection all closely approximate the Gauss-Legendre points! This figure---and the numerous other cases we tested for $d=1$---show that the difference between the various optimization strategies is not incredibly significant; all of them tend to converge close to the optimal  solution. 

But what happens when $d$ increases? Figure \ref{gauss2d}(a) plots the subsampled points obtained from the three linear algebra optimization strategies when subsampling from a tensor grid with order 50 in both $\zeta_1$ and $\zeta_2$ directions. The basis used---i.e., the columns of $\mA$---is an isotropic tensor product with an order of 3. As expected, the obtained points are very similar to Gauss-Legendre tensor product points with order 3. In Figure \ref{gauss2d}(b), instead of subsampling from a higher order tensor grid, we subsample from Christoffel samples, which in the uniform case corresponds to the Chebyshev distribution; the basis is the same as in (a). In this case it is interesting to note that the LU and QR pivoting approaches, which are similar, yield slightly different results when compared to the SVD approach.
\end{svgraybox}

The use of such rank revealing factorizations for minimizing the condition number of a Vandermonde-type matrix have been previously investigated in \cite{bos2010computing, seshadri2017effectively} albeit on different meshes. When pruning rows from $\mA$, the condition number of $\hat{\mA}_{\vz}$ can be bounded by
\begin{equation}
\kappa\left(\hat{\mA}_{\vz} \right)\leq\kappa\left(\mA\right)\sqrt{1+s^{2}n\left(m-n\right)},
\end{equation}
where $s > 1$ is a constant based on the pivoting algorithm. In Seshadri et al. \cite{seshadri2017effectively} the authors draw upon a subset selection strategy from Chapter 5 of \cite{golub2012matrix} that uses both the SVD and a rank revealing QR factorization for identifying a more well conditioned submatrix. However, as $m$ and $n$ get progressively larger, the costs associated with computing both the SVD $\mathcal{O}(m^2n^3)$ and a Householder (or modified Gram-Schmidt)-based rank revealing QR factorization should be monitored. Other rank revealing techniques such as strong LU factorizations \cite{miranian2003strong} and randomized QR with column pivoting \cite{martinsson2017householder, duersch2017randomized} may also be opted for. 

To speed up the optimization process---a consideration that may be ignored when model evaluations themselves take hours---other techniques have been proposed. Shin and Xiu \cite{shin2016nonadaptive, shin2016near} optimize a metric similar to P4
\begin{equation}
\underset{\vz}{\text{maximize}}\left(\frac{\sqrt{\text{det}\left(\hat{\mA}_{\vz}^{T}\hat{\mA}_{\vz}\right)}}{\prod_{i=1}^{k}\left\Vert \hat{\mA}_{\vz}(i)\right\Vert }\right)^{1/k},
\end{equation}
where $\hat{\mA}(i)_{\vz}$ denotes the $i-$th column of $\hat{\mA}_{\vz}$. To ascertain which rows $\vz$ to use, the authors outline a greedy strategy that begins with a few initially selected rows and then adds rows based on whether they increase the determinant of the Gramian, $\hat{\mA}_{\vz}^{T}\hat{\mA}_{\vz}$. To abate the computational cost of computing the determinant for each candidate, they devise a rank-2 update to the determinant (of the Gramian) based on Sylvester's determinant formula and the Sherman-Morrison formula. In a similar vein, Ghili and Iaccarino \cite{ghili2017least} lay out a series of arguments---motivated by reducing the operational complexity of the determinant-based optimization---for optimizing the trace of the design matrix instead. They observe that since the Frobenius norm is an upper bound on the spectral norm and because all the singular values contribute to the aliasing error, optimizing over this norm will undoubtedly yield a small spectral norm. The authors formulate a coordinate descent optimization approach to find suitable points and weights. 

\subsection{Convex relaxation via Newton's method}
\label{sec:convex}
In general the aforementioned objectives are non-convex. We discuss an idea that permits these objectives to be recast as a convex optimization problem. Our idea originates from the sensor selection problem \cite{hovland1997dynamic, kincaid2002d}: given a collection of $m$ sensor measurements---where each measurement is a vector---select $k$ measurements such that the error covariance of the resulting ellipsoid is minimized. We remark here that while a generalized variant of the sensor selection problem is NP-hard (see \cite{bian2006utility}); the one we describe below has not yet been proven to be NP-hard. Furthermore, by selecting the optimization variable to be a boolean vector $\vz$ that restricts the measurements selected, this problem can be cast as a determinant maximization problem \cite{vandenberghe1998determinant} where the objective is a concave function in $\vz$ with binary constraints on its entries. This problem can be solved via interior point methods and has complexity $\mathcal{O}(m^3)$. Joshi and Boyd \cite{joshi2009sensor} provide a formulation for a \emph{relaxed sensor selection problem} that can be readily solved by Newton's method, where the binary constraint can be substituted with a penalty term added to the objective function. By replacing their sensor measurements with rows from our Vandermonde type matrix, one arrives at the following maximum volume problem
\begin{align}
\begin{split}
\underset{\vz \in \mathbb{R}^{M}}{\text{minimize}}  \; \; \; &  -\text{log} \; \text{det}\left( \frac{1}{\tau} \sum_{i=1}^{M}  z_{i} \va_{i}^{T} \ve  \right) -\lambda\sum_{i=1}^{M}\left( \text{log} \left(z_{i}\right)+ \text{log} \left(1-z_{i}\right)\right)\\
\text{subject to}  \; \; \; & \mathbf{1}^{T}\vz =K \\
 & 0 \leq z_i \leq 1, \; i=1, \ldots, M.
 \end{split}
 \label{cvx}
\end{align}
where the positive constant $\lambda$ is used to control the quality of the approximation. Newton's approach has complexity $\mathcal{O}(m^3)$ with the greatest cost arising from the Cholesky factorization when computing the inverse of the Hessian \cite{joshi2009sensor}. The linear constraint is solved using standard KKT conditions as detailed in page 525 of \cite{boyd2004convex}. In practice this algorithm requires roughly 10-20 iterations and yields surprisingly good results for finding suitable quadrature points.

\begin{figure}
\begin{subfigmatrix}{2}
\subfigure[]{\includegraphics[]{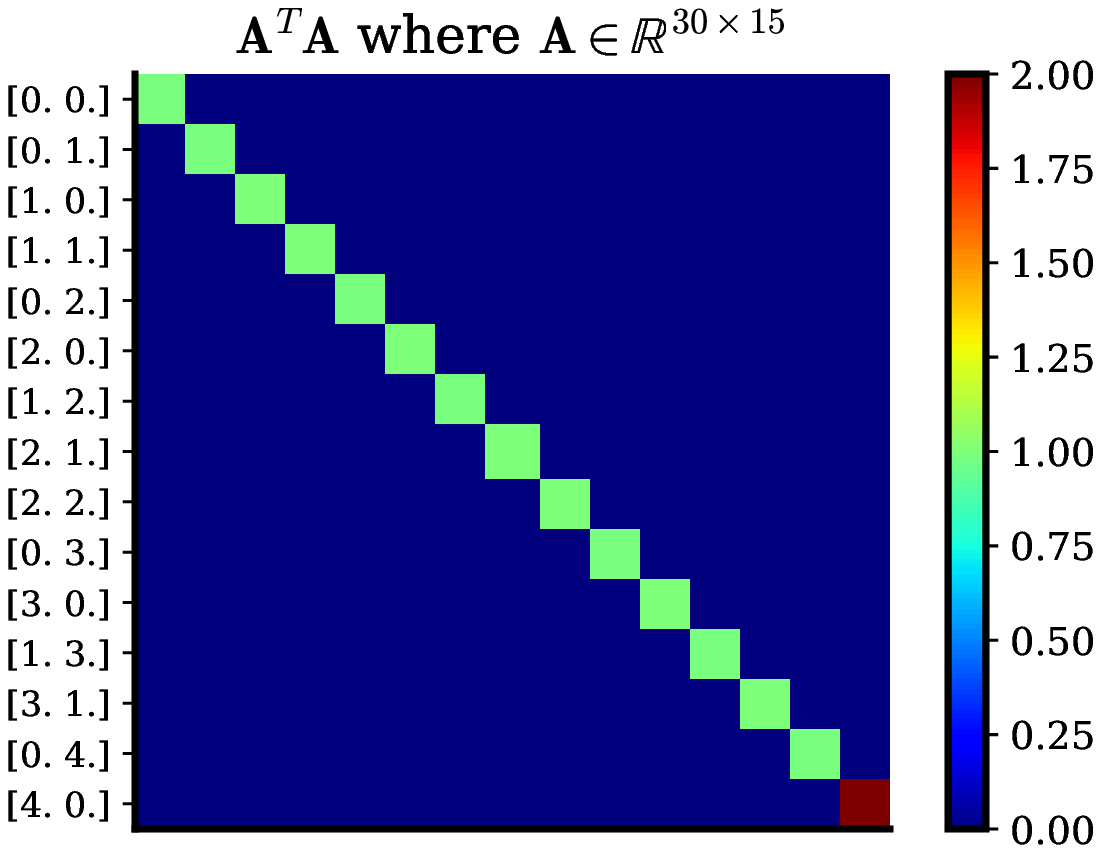}}
\subfigure[]{\includegraphics[]{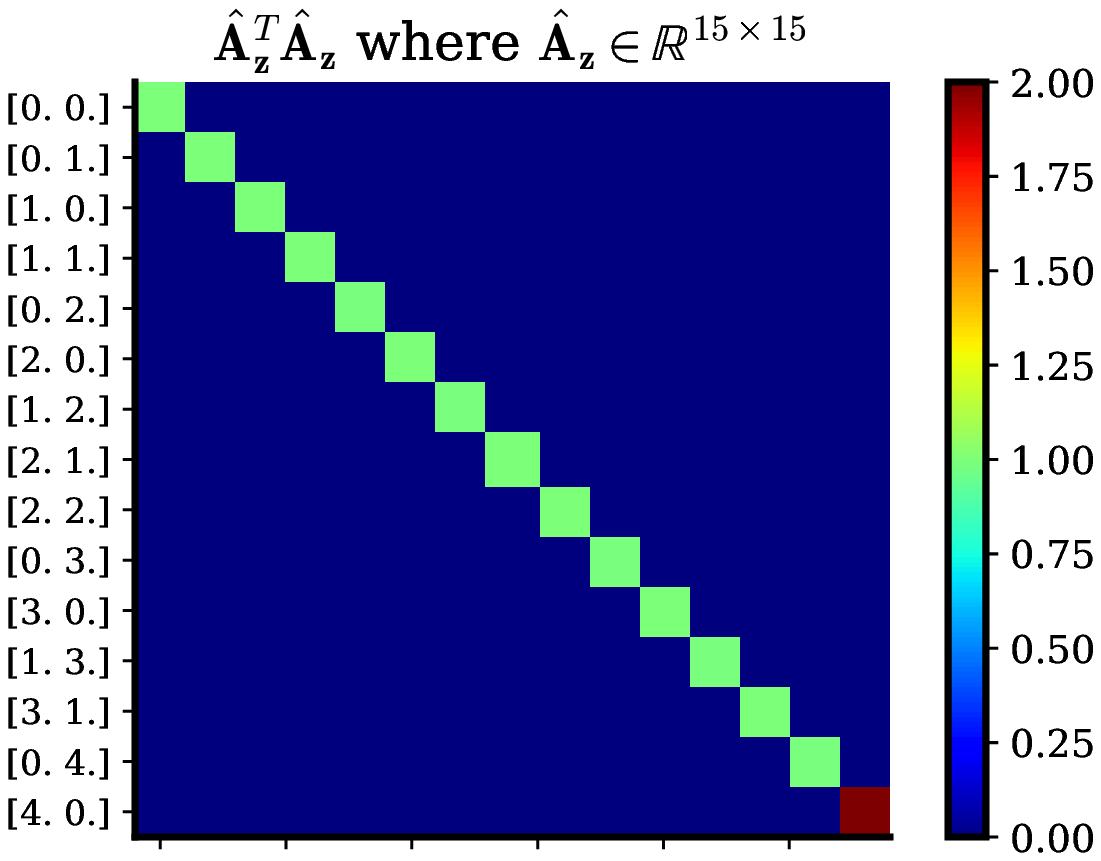}}
\subfigure[]{\includegraphics[]{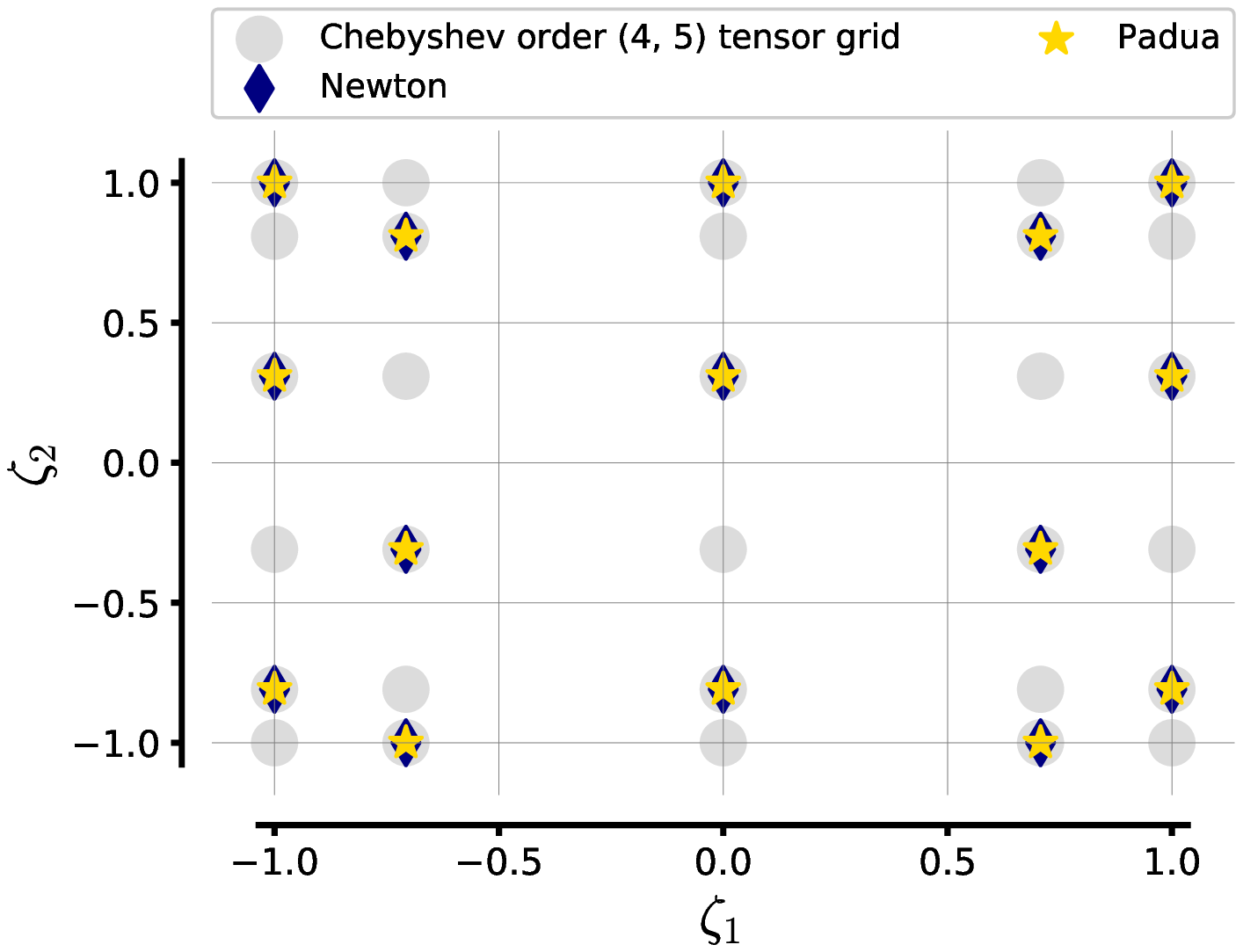}}
\end{subfigmatrix}
\caption{The Gram matrix associated with Chebyshev tensor product grid with $4\times 5=30$ points and a total order basis with maximum degree of 4 in (a). The subsampled Gram matrix via convex optimization in (b) and a comparison of the samples and subsamples in (c).}
\label{padua}
\end{figure}

\begin{svgraybox}
\textbf{Padua points via Convex optimization:}
It is difficult to select points suitable for interpolation with total order polynomials, that are \emph{unisolvent}---a condition where the points guarantee a unique interpolant \cite{trefethen2017cubature}. One such group of points that does admit unisolvency are the famous \emph{Padua points}. For a positive integer $N$, these points are given by $(\zeta^{(1)}_m, \zeta^{(2)}_m) \in [-1,1]^2$ where
\begin{equation}
\zeta^{(1)}_{m}=cos\left(\frac{\pi\left(m-1\right)}{N}\right), \zeta^{(2)}_{m}=\begin{cases}
\begin{array}{c}
cos\left(\frac{\pi\left(2k-1\right)}{N-1}\right)\\
cos\left(\frac{\pi\left(2k-2\right)}{N-1}\right)
\end{array} & \begin{array}{c}
m\; \textrm{odd}\\
m\; \textrm{even}
\end{array}
\end{cases}
\end{equation}
and $1\leq m \leq N+1$, $1\leq k \leq 1 + N/2$ \cite{caliari2011padua2dm}. The Padua points have a provably minimal growth rate of $\mathcal{O}(log ^2 N)$ of the Lebesgue constant, far lower than Morrow-Patterson or the Extended Morrow-Patterson points \cite{bos2006bivariateb}. Two other characterizations may also be used to determine these points; they are formed by the intersection of certain Lissajous curves and the boundary of $[-1, 1]^2$ or alternatively, every other point from an $(N+1) \times (N+2)$ tensor product Chebyshev grid \cite{bos2006bivariate, bos2007bivariate}. 

As an example, consider the case where $N=4$ resulting in a 30-point Chebyshev tensor grid and a total order basis where the highest order is 4. Figure~\ref{padua}(a) plots the corresponding gram matrix $\mG = \mA^T \mA$. As expected the quadrature rule can integrate all polynomials except for the one corresponding to the last row, i.e., $\bm{\psi_{(4,0)}}$, whose integrand has an order of 8 along the first direction $\zeta_{1}$, beyond the degree of exactness of the 5 points. What is fascinating is that when these points are subsampled one can get a subsampled Gram matrix i.e., $\hat{\mA}_{\vz}^{T} \hat{\mA}_{\vz}$ to have the same degree of exactness, as illustrated in Figure~\ref{padua}(b). To obtain this result, the aforementioned convex optimization via Newton's method was used on $\mA$, yielding the Padua points---i.e., every alternate point from the 30-point Chebyshev tensor grid. 
\end{svgraybox}

\begin{figure}
\begin{center}
\includegraphics[scale=0.5]{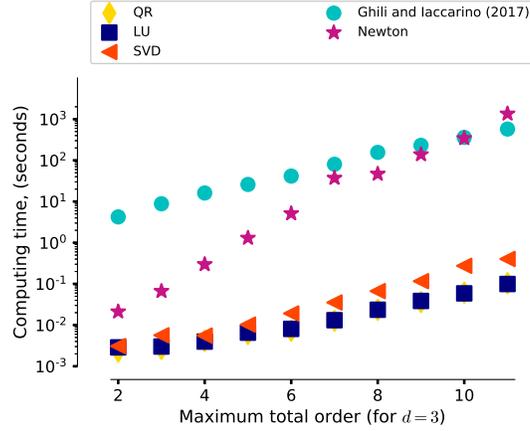}
\caption{Representative computing times for the different optimization strategies for $m=n$ and $d=3$ using Chebyshev samples and a total order basis.}
\label{experiment3D_analysis}
\end{center}
\end{figure}

\begin{svgraybox}
\textbf{Comparing execution time:}
In our numerical experiments on some of the aforementioned optimization strategies, we have alluded to the operational complexity. We make an obvious remark here, that as the both the polynomial degree and the dimension increases, the computing time rises. Non-linear algebra approaches such as \cite{ghili2017least} and the Newton approach presented here, have other tunable parameters that can either add to, or decrease computational run time.

Figure \ref{experiment3D_analysis} plots the running times for $d=3$ and $m=n$ when subsampling from grid of Chebyshev points---where the number of points is given by the $\left( \text{maximum total order} + 1\right)^{d}$---using the different optimization strategies. These experiments were carried out on a 3.1 GHz i7 Macbook Pro with 16GB of RAM. For completeness, a summary of the operational complexity of these techniques is provided in \ref{operational_complexity}. The condition numbers of the matrices obtained from QR, LU and the SVD approach were generally lower (similar order of magnitude) than those obtained from the other two techniques. 
\end{svgraybox}

\begin{table}
\begin{center}
\caption{Summary of the operational complexity of some of the optimization strategies.}
\vspace{2 mm}
\begin{tabular}{l l }
\hline
\textbf{Optimization strategy} & \textbf{Complexity} \\
\hline
SVD-based subset selection & $\mathcal{O} \left( nm^2  \right) + \mathcal{O} \left( nm^2  \right)$   \\
QR column pivoting & $ \mathcal{O} \left( nm^2  \right)$  \\
LU row pivoting & $ \mathcal{O} \left( m^2  \right)$ \\
Newton's method & $ \mathcal{O} \left( n^3  \right)$ \\
Frobenius norm optimization (Ghili and Iaccarino 2017) & $ \mathcal{O} \left( dnm^2  \right)$ \\          
\hline
\end{tabular}
\end{center}
\label{operational_complexity}
\end{table}

\subsection{Moment-based optimization strategies}
To motivate this section, consider the integration property of orthogonal polynomials, highlighted in \eqref{doe_1}. In the univariate case, one can write 
\begin{align}
\begin{split}
& \mP \vw = \ve \\
\text{where} \; \; \;  & \mP(i, j) = \psi_{i} \left( \zeta_{j} \right), \; \; \; \vw(i) = \omega_{i}, \; \; \; \text{and} \; \; \;  \ve^{T} = \left[ 1, 0, \ldots, 0 \right].
\end{split}
\end{align}
Assuming $\mP \in \mathbb{R}^{m \times m}$ where one uses the first m Gauss-Legendre  points in $\zeta_{j}$ and the first $(m-1)$ order Legendre polynomials in $\psi_{i}$, then the solution $\vw$ for the above linear set of equations will yield the first $m$ Gauss-Legendre weights. In the case where the points are not from a known quadrature rule---i.e., if they are from a random distribution---then one has to ensure that the weights are positive via a constraint. 

Instead of working with matrices and optimizing for the objectives in Table~\ref{table1}, one can frame an optimization problem centered on the computation of moments where one optimizes over the space of quadrature weights (non-negative). This \emph{moment-based approach} for finding quadrature points and weights is utilized in Keshavarzzadeh et al. \cite{keshavarzzadeh2018numerical} where the authors solve a \emph{relaxed} version of this constraint linear system by minimizing $\left\Vert \mathbf{P}\mathbf{w}-\mathbf{e}\right\Vert_{2}$. In Ryu and Boyd \cite{ryu2015extensions}, the authors present numerical quadrature as a solution to the infinite-dimensional linear program
\begin{align}
\begin{split}
\underset{\vv \in \mathbb{R}^{m}}{\text{minimize}}  \; \; \; &  \vf^{T} \vv \\
\text{subject to}  \; \; \; &  \vp_{j}^{T} \vv = c_{i}, \; \;  \; \text{where} \; \; \;    i=1, \ldots, n.
 \end{split}
 \label{lp}
\end{align}
Here components of $\vf \in \mathbb{R}^{m}$ are given by $\vf \left( i \right) = f \left(\bm{ \zeta}_i \right)$, and $\vp_{j} \in \mathbb{R}^{m}$ has components $\vp_{j} \left( i \right) = \bm{\psi}_{j} \left( \bm{\zeta}_{i} \right)$; the optimization variable $\vv \left( i \right) = \omega_{i}$ represents the quadrature weights. The constants $c_i$ in the equality constraints are determined analytically, as they involve integrating a known polynomial over the support of $f$. The problem can be interpreted as a weighted $l_1$ optimization problem, as we require $\vv$ to have as many zeros as possible and yet satisfy the above constraints. As this problem is NP-hard, Ryu and Boyd propose a two-stage approach to solve it; one for generating an initial condition and another for optimizing over $\vv$. Their approach has been recently adapted in Jakeman and Narayan \cite{jakeman2017generation} who propose least absolute shrinkage and selection operator (LASSO) for finding the initial condition. They then proceed to solve the optimization using a gradient-based nonlinear least squares optimizer. Their results are extremely promising---numerical experiments using their technique show orders of magnitude improvement in convergence compared to tensor and sparse grid rules.

\section{Concluding remarks and the future}
In this paper, we provided an overview of strategies for finding quadrature points using ideas from polynomial least squares. Although we have sought to keep our review as detailed as possible, readers will forgive us for omitting various techniques. For example, we did not discuss optimal design of experiment based samples, although we point the interested reader to the review offered in \cite{hadigol2017least}; for convex relaxations of the various optimal design of experiment problems we refer the reader to section 7.5.2 of Boyd and Vandenberghe \cite{boyd2004convex}. In addition, we have also highlighted a new convex optimization strategy that uses Newton's method for finding the best subsamples by maximizing the volume of the confidence ellipsoid associated with the Vandermonde-type matrix. 

So what’s next? Based on our review we offer a glimpse of potential future areas of research that could prove to be fruitful:
\begin{enumerate}
\item Randomized greedy linear algebra approaches for finding suitable quadrature samples. Existing approaches are tailored for finding pivot columns for tall matrices; for our problem we require these techniques to be applicable to fat matrices.
\item Large scale (and distributed) variants of the convex optimization strategies detailed, including an alternating direction method of multiplies (ADMM) \cite{boyd2011distributed} formulation for the Newton’s method technique presented in this paper;
\item Heuristics for optimizing the weights when the joint density of the samples is not known---a problem that arises in data science; typically in uncertainty quantification the joint density $\bm{\rho}$ is assumed to be known.
\item The development of an open-source repository of near-\emph{optimal} points for the case where $m=n$ for different total order basis and across different $d$. 
\item Building on (1) and following the recent work by Shin and Xiu \cite{shin2017randomized} and Wu et al. \cite{wu2017randomized}, the use of approaches such as the randomized Kaczmarz algorithm \cite{strohmer2009randomized} for solving the least squares problem in \eqref{leastsquares}. The essence of the idea here is that as $d$ and the highest multivariate polynomial degree get larger, the matrix $\hat{\mathbf{A}}_{\vz}$ can not be stored in memory---a requirement for standard linear algebra approaches. Thus techniques such as the Kaczmarz algorithm, which solve for $\vx$ by iteratively requiring access to rows of $\hat{\mathbf{A}}_{\vz}$ and elements of $\vb$, are useful.

\end{enumerate}

\section*{Acknowledgements}
This work was carried out while PS was visiting the Department of Mechanical, Chemical and Materials Engineering at Universit\'{a} di Cagliari in Cagliari, Sardinia; the financial support of the University's Visiting Professor Program is gratefully acknowledged. The authors are also grateful to Akil Narayan for numerous discussions on polynomial approximations and quadratures. 

\bibliography{references}
\bibliographystyle{plain}
\end{document}